\documentclass[12pt]{smfart}

\usepackage[latin1]{inputenc}
\usepackage[francais]{babel}
\usepackage[leqno]{amsmath}
\usepackage{smfthm}
\usepackage{amsfonts}
\usepackage{amssymb}
\usepackage{eucal}
\usepackage{mathrsfs}

\def\C{\mathbb{C}}
\def\R{\mathbb{R}}
\def\N{\mathbb{N}}
\def\Q{{\mathbb{Q}}_p}
\def\Qpbar{{\overline{\mathbb{Q}}}_p}
\def\Z{{\mathbb{Z}}_p}

\def\O{\mathcal{O}_L}
\def\r{\mathscr{R}_L}

\def\G{{\rm GL}_2(\Q)}
\def\B{{\rm B}(\Q)}
\def\K{{\rm GL}_2(\Z)}
\def\BK{{\rm B}(\Z)}

\def\={\buildrel{\rm d\acute ef}\over =}

\def\g{{\rm Gal}(\Qpbar/\Q)}
 
\def\atplus{\widetilde{\mathbb{A}}^+} 
\def\bcris{\mathbb{B}_{\rm cris}}

\def\bdR{\mathbb{B}_{\rm dR}} 
\def\dcris{D_{\rm cris}} 
\def\ddR{D_{\rm dR}} 
\def\calE{\mathscr{E}_L}
\def\calO{\mathscr{O}_{\mathscr{E},L}}
\def\calR{\mathscr{R}_L}

\author[L. Berger]{Laurent Berger}
\address{C.N.R.S. \& I.H.\'E.S.\\
Le Bois-Marie\\
35 route de Chartres\\
91440 Bures-sur-Yvette \\ 
France}
\email{lberger@ihes.fr}
\urladdr{www.ihes.fr/\~{}lberger/}

\author[C. Breuil]{Christophe Breuil}
\address{C.N.R.S. \& I.H.\'E.S.\\
Le Bois-Marie\\
35 route de Chartres\\
91440 Bures-sur-Yvette \\ 
France}
\email{breuil@ihes.fr}
\urladdr{www.ihes.fr/\~{}breuil/}

\title[Repr\'esentations cristallines irr\'eductibles de ${\rm
    GL}_2(\Q)$]{Repr\'esentations cristallines irr\'eductibles de
  ${\rm GL}_2(\Q)$}  

\alttitle{Irreducible crystalline representations of 
${\rm GL}_2(\Q)$}

\subjclass{11F}

\NumberTheoremsIn{subsection}

\begin{document}

\begin{abstract}
Dans \cite[\S1.3]{Br2}, des 
repr\'esentations unitaires de $\G$ sur des espaces de Banach 
$p$-adiques sont associ\'ees aux repr\'esentations 
cristallines irr\'educti\-bles de dimension $2$ de $\g$. 
Certaines conjectures sur ces Banach y sont formul\'ees 
(non nullit\'e, irr\'eductibilit\'e topologique, admissibilit\'e). Nous 
d\'emontrons ces conjectures en r\'einterpr\'etant 
ces espaces de Banach comme 
espaces de fonctions sur $\Q$ d'un certain type, puis en utilisant la 
th\'eorie des $(\varphi,\Gamma)$-modules associ\'es aux repr\'esentations 
cristallines correspondantes.
\end{abstract}

\begin{altabstract}
In \cite[\S1.3]{Br2}, 
some unitary representations of $\G$ on $p$-adic Banach spaces are
associated to $2$-dimensional irreducible crystalline representations
of $\g$. Some conjectures are formulated concerning those Banach
spaces (non triviality, topological irreducibility, admissibility). We
prove those conjectures by reinterpreting those Banach as spaces of
functions of a certain type on $\Q$, and then by using the theory of 
$(\varphi,\Gamma)$-modules associated to the corresponding crystalline
representations. 
\end{altabstract}

\maketitle

\setcounter{tocdepth}{2}
\tableofcontents

\setlength{\baselineskip}{18pt}

\section{Introduction et notations}

\subsection{Introduction}

Soit $p$ un nombre premier et $n\in \N$. Dans la recherche d'une 
correspondance \'eventuelle entre (certaines) repr\'esentations 
$p$-adiques $V$ de $\g$ de dimension $n$ et (certaines) 
repr\'esentations $p$-adiques $\Pi(V)$ de ${\rm GL}_n(\Q)$, un des 
premiers cas \`a regarder est certainement celui o\`u $n=2$ et $V$ est 
cristalline, absolument irr\'eductible et $\varphi$-semi-simple. Cette 
derni\`ere condition signifie que le Frobenius $\varphi$ sur le
$\varphi$-module  
filtr\'e $\dcris(V)$ associ\'e par Fontaine \`a $V$ est semi-simple.
La repr\'esentation $V$ a des poids de Hodge-Tate distincts $i_1< i_2$ et, si 
l'on note ${\rm Alg}(V)$ la repr\'esentation alg\'ebrique de $\G$ de 
plus haut poids $(i_1,i_2-1)$ et ${\rm Lisse}(V)$ la repr\'esentation 
lisse irr\'eductible de $\G$ associ\'ee par la correspondance locale de 
Hecke \`a la repr\'esentation non ramifi\'ee de $\g$ envoyant un 
Frobenius arithm\'etique sur $\varphi^{-1}$, la repr\'esentation 
$\Pi(V)$ est simplement le compl\'et\'e $p$-adique de la 
repr\'esentation localement alg\'ebrique ${\rm Alg}(V)\otimes {\rm 
Lisse}(V)$ par rapport \`a un r\'eseau stable par $\G$ et de type fini 
sous l'action de $\G$. Ainsi, $\Pi(V)$ est un espace de Banach 
$p$-adique muni d'une action continue unitaire de $\G$ (i.e. laissant 
une norme invariante). Notons que, lorsque ${\rm Lisse}(V)$ est de 
dimension $1$, cette d\'efinition de $\Pi(V)$ doit \^etre modifi\'ee.

Le probl\`eme est qu'il n'est pas du tout \'evident qu'un tel r\'eseau 
existe, ou, de mani\`ere \'equivalente, que $\Pi(V)$ soit non nul ! 
Dans \cite[th\'eor\`eme 1.3]{Br1}, la non nullit\'e de $\Pi(V)$ est 
d\'emontr\'ee pour $i_2-i_1<2p$ (essentiellement), et dans
\cite[th\'eor\`eme 1.3.3]{Br2}, 
son admissibilit\'e (au sens de \cite{ST3}) et son irr\'eductibilit\'e 
topologique (avec une condition suppl\'ementaire pour cette 
derni\`ere). La m\'ethode repose sur le calcul de la r\'eduction d'une 
boule unit\'e de $\Pi(V)$ modulo $p$ (plus exactement modulo l'id\'eal 
maximal des coefficients). La restriction ci-dessus sur les poids n'est 
pas th\'eorique mais simplement due au fait que les calculs deviennent 
de plus en plus difficiles quand les poids grandissent.

Lorsque $V$ est de dimension $2$, absolument irr\'eductible mais cette 
fois semi-stable non-cristalline, $\Pi(V)$ est d\'efini dans \cite{Br2} 
et \cite{Br3} et des conjectures analogues (non-nullit\'e, 
admissibilit\'e, etc.) formul\'ees et tr\`es partiellement 
d\'emontr\'ees via la th\'eorie globale dans \cite{Br3} ou via des calculs 
de r\'eduction modulo $p$ sur le demi-plan $p$-adique dans
\cite{BM}. P. Colmez dans \cite{Co2} a vu que la th\'eorie des  
$(\varphi,\Gamma)$-modules de Fontaine permettait de d\'emontrer 
\'el\'egamment ces conjectures dans le cas semi-stable en construisant 
un mod\`ele de la restriction au Borel sup\'erieur de la 
repr\'esentation duale $\Pi(V)^*$ \`a partir du 
$(\varphi,\Gamma)$-module de $V$.

Il \'etait donc naturel de regarder si un tel mod\`ele existait aussi 
dans le cas {\it a priori} plus simple des repr\'esentations $V$ 
cristallines et s'il permettait de d\'emontrer les conjectures 
ci-dessus. La r\'eponse est affirmative et fait l'objet du pr\'esent 
article. Comme dans \cite{Co2}, on d\'emontre donc un isomorphisme 
Borel-\'equivariant entre $\Pi(V)^*$ et  $(\varprojlim_{\psi} 
D(V))^{\rm b}$ (th\'eor\`eme \ref{phigammabanach}) o\`u $D(V)$ est le 
$(\varphi,\Gamma)$-module associ\'e \`a la repr\'esentation
cristalline $V$ et o\`u la limite projective consiste en les suites
$\psi$-compatibles born\'ees  
d'\'el\'ements de $D(V)$. Pour montrer cet isomorphisme, il est 
n\'ecessaire de passer par une autre description interm\'ediaire de 
$\Pi(V)$ comme espace fonctionnel $p$-adique, i.e. comme un espace de 
fonctions continues sur $\Q$ d'un certain type (th\'eor\`eme 
\ref{complete}). Les r\'esultats c\^ot\'e $(\varphi,\Gamma)$-modules 
permettent alors de d\'eduire le r\'esultat principal de cet article 
(\S\ref{resul}) :

\begin{enonce*}{Th\'eor\`eme}
Supposons $V$ de dimension $2$, cristalline, absolument 
irr\'e\-ductible et $\varphi$-semi-simple, alors $\Pi(V)$ est non nul, 
topologiquement irr\'eductible et admissible.
\end{enonce*}

On obtient aussi deux autres corollaires, l'un concernant tous les 
r\'eseaux possibles stables par $\G$ dans ${\rm Alg}(V)\otimes {\rm 
Lisse}(V)$ (corollaire \ref{reseaurigolo}), l'autre concernant les 
vecteurs localement analytiques dans $\Pi(V)$ (corollaire \ref{anal}).

Une version pr\'eliminaire des r\'esultats de cet article
a fait l'objet d'un cours \`a Hangzhou en ao\^ut 2004 et 
les notes de ce cours sont disponibles (en anglais, voir \cite{BB}).  

Terminons cette introduction en insistant sur le fait que des 
\'enonc\'es non triviaux sur des repr\'esentations de $\G$ 
(irr\'eductibilit\'e des $\Pi(V)$ etc.) se d\'emontrent finalement en 
passant par des repr\'esentations de $\g$.

\subsection{Notations}

On fixe $\Qpbar$ une cl\^oture alg\'ebrique de $\Q$, 
{\og ${\rm val}$ \fg} la valuation sur $\Qpbar$ 
telle que ${\rm val}(p)\=1$, $|\cdot |$ la norme $p$-adique
$|x|\=p^{-{\rm val}(x)}$ et $\C_p$ le compl\'et\'e de
$\Qpbar$ pour $|\cdot|$. On normalise l'isomorphisme de la th\'eorie
du corps de classes local en envoyant les uniformisantes sur les
Frobenius g\'eom\'etriques. On note $\varepsilon$ le caract\`ere
cyclotomique $p$-adique de $\g$ vu aussi comme caract\`ere de
$\Q^{\times}$. En particulier, $\varepsilon(p)=1$ et
$\varepsilon\!\!\mid_{\Z^{\times}}:\Z^{\times}\rightarrow \Z^{\times}$
est l'identit\'e. On note ${\rm nr}(x)$ le caract\`ere non ramifi\'e
de $\Q^{\times}$ envoyant $p$ sur $x$. On d\'esigne par $L$ une
extension finie de $\Q$, $\O$ son anneau d'entiers, $k_L$ son corps
r\'esiduel et $\pi_L$ une uniformisante. On note $\Q(\pmb{\mu}_{p^n})$
l'extension finie de $\Q$ dans $\Qpbar$ engendr\'ee par les racines
$p^n$-i\`emes de l'unit\'e et
$\Q(\pmb{\mu}_{p^{\infty}})\=\cup_{n \geq 0} 
\Q(\pmb{\mu}_{p^n})$. On note $V$
une repr\'esentation $p$-adique de $\g$, c'est-\`a-dire un $L$-espace
vectoriel de dimension finie muni d'une action lin\'eaire et continue
de $\g$, et $T$ un $\O$-r\'eseau de $V$ stable par $\g$. Enfin, $\B$
d\'esigne les matrices triangulaires sup\'erieures dans $\G$ et $\BK$
le sous-groupe de ces matrices qui sont dans $\K$.

\section{S\'eries formelles et $(\varphi,\Gamma)$-modules}\label{phigamma}

Le but de cette partie est de donner des pr\'eliminaires sur les
$(\varphi,\Gamma)$-modules $D(V)$ et leurs treillis $D^{\sharp}(V)$ qui
sont utilis\'es dans les parties suivantes 
pour l'\'etude des repr\'esenta\-tions
cristallines de $\G$.  

\subsection{Quelques anneaux de s\'eries formelles}\label{series}

Le but de ce paragraphe est d'introduire certains anneaux de
s\'eries formelles ($\calO$, $\calE$ et $\calR^+$), ainsi que
certaines des structures dont ils sont munis et dont nous avons besoin
dans la suite de cet article. 

\begin{itemize}
\item[(i)] On note $\calO$ l'anneau form\'e des
s\'eries $\sum_{i \in \mathbb{Z}} a_i X^i$ telles que $a_i \in \O$ et
$a_{-i} \rightarrow 0$ quand $i \rightarrow +\infty$. C'est un anneau
local de corps r\'esiduel $k_L((X))$.
\item[(ii)] On note $\calE \= \calO[1/p]$, c'est un
corps local de dimension $2$.
\item[(iii)] On note $\calR^+$ l'anneau des 
s\'eries formelles $f(X) \in
L[[X]]$ qui convergent sur le disque unit\'e, ce qui fait que
$L \otimes_{\O} \O[[X]]$ s'identifie au sous-anneau 
de $\calR^+$ constitu\'e des
s\'eries formelles \`a coefficients born\'es.
\end{itemize}

Le corps $\calE$ est muni de la norme de Gauss $\|\cdot\|_{\rm Gauss}$
d\'efinie par $\|f(X)\|_{\rm Gauss} \= \sup |a_i|$ si $f(X)=\sum_{i
  \in \mathbb{Z}} a_i X^i$. L'anneau des entiers de $\calE$ pour cette
norme est $\calO$, et la norme de Gauss induit sur $\calO$ la
topologie $\pi_L$-adique. L'application naturelle $\calO \rightarrow
k_L((X))$ est alors continue, si l'on donne \`a $\calO$ la topologie
$\pi_L$-adique et \`a $k_L((X))$ la topologie discr\`ete. 

On peut d\'efinir une topologie moins fine sur $\calO$, la topologie
faible, d\'efinie par le fait que les $\{\pi_L^i \calO + X^j
\O[[X]]\}_{i,j \geq 0}$ forment une base de voisinages
de z\'ero, et la topologie faible sur $\calE = \cup_{k \geq 0}
\pi_L^{-k} \calO$ qui est la topologie de la limite inductive. Cette
topologie induit la topologie $(\pi_L,X)$-adique sur $\O[[X]]$, et 
l'application naturelle $\calO \rightarrow
k_L((X))$ est alors continue, si l'on donne \`a $\calO$ la topologie
faible et \`a $k_L((X))$ la topologie $X$-adique. 

Si $\rho < 1$, alors on peut d\'efinir une norme $\|\cdot\|_{D(0,\rho)}$
sur $\calR^+$ par la formule :
\begin{equation}\label{normerho}
\|f(X)\|_{D(0,\rho)} \= \sup_{\substack{z \in \C_p \\ |z| \leq \rho}}
|f(z)|  =  \sup_{i \geq 0} |a_i| \rho^i,  
\end{equation}
si $f(X)=\sum_{i \geq 0} a_i X^i$. L'ensemble des normes $\{
\|\cdot\|_{D(0,\rho)} \}_{0 \leq \rho < 1}$ d\'efinit une topologie sur
$\calR^+$ qui en fait un espace de Fr\'echet.

\begin{defi}\label{dordrer}
Si $f(X) \in \calR^+$ et $r \in \R_{\geq 0}$, on dit que $f(X)$ est d'ordre $r$ (il serait plus correct de dire {\og d'ordre $\leq r$ \fg}) si pour un $\rho$ tel que $0 < \rho < 1$, la suite $\{ p^{-nr} \|f(X)\|_{D(0,\rho^{1/p^n})} \}_{n \geq 0}$ est
born\'ee.
\end{defi}

Il est facile de voir que si c'est vrai pour un choix de $0
< \rho < 1$, alors c'est vrai pour tout choix de $0 < \rho < 1$. 
Un exemple de s\'erie d'ordre $1$ est donn\'e par $f(X)=\log(1+X)$.

Soit $\Gamma \= {\rm Gal}(\Q(\pmb{\mu}_{p^\infty})/\Q)$. Les anneaux
$\calO$, $\calE$ et $\calR^+$ d\'efinis
ci-dessus sont munis d'une action de $\Gamma$, telle que si $\gamma
\in \Gamma$, alors $\gamma$ agit par un morphisme de $L$-alg\`ebres,
et $\gamma(X) \= (1+X)^{\varepsilon(\gamma)}-1$ o\`u $\varepsilon :
\Gamma \rightarrow \Z^{\times}$ est le caract\`ere cyclotomique. 
On v\'erifie facilement que $\Gamma$ agit par des isom\'etries, pour
toutes les normes et topologies d\'efinies ci-dessus. Ces anneaux
sont aussi munis d'un morphisme de Frobenius $\varphi$, qui est
lui-aussi un morphisme de $L$-alg\`ebres, tel que $\varphi(X) \=
(1+X)^p-1$. Cette application est continue pour toutes les topologies
ci-dessus et commute \`a l'action de $\Gamma$. Remarquons enfin que
$\varphi$ et l'action de $\Gamma$ stabilisent $L \otimes_{\O}
\O[[X]]$. 

L'anneau $\calO$ est un $\varphi(\calO)$-module libre de rang $p$,
dont une base est donn\'ee par $\{(1+X)^i\}_{0 \leq i \leq p-1}$. Si
$x \in \calO$, on peut donc \'ecrire $x=\sum_{i=0}^{p-1} (1+X)^i
\varphi(x_i)$.

\begin{defi}\label{defpsianno}
On d\'efinit un op\'erateur $\psi : \calO \rightarrow
\calO$ par la formule $\psi(x) \= x_0$ si $x=\sum_{i=0}^{p-1} (1+X)^i
\varphi(x_i)$. On d\'efinit de mani\`ere similaire un
op\'erateur $\psi : \calR^+ \rightarrow \calR^+$. 
\end{defi}

Cet op\'erateur v\'erifie alors $\psi(\varphi(x)y) = x \psi(y)$ et
commute \`a l'action de $\Gamma$. Il ne commute pas \`a $\varphi$ et
n'est pas $\calO$-lin\'eaire. 

\subsection{Repr\'esentations $p$-adiques et
$(\varphi,\Gamma)$-modules} \label{pgmod}

Rappelons qu'une repr\'esen\-tation $L$-lin\'eaire est un $L$-espace
vectoriel $V$ de 
dimension finie muni d'une action lin\'eaire et continue de $\g$. Il
est assez difficile de d\'ecrire un tel objet, et nous allons voir
dans ce paragraphe qu'une certaine classe de {\og
$(\varphi,\Gamma)$-modules sur $\calE$ \fg} permet d'en donner une
description explicite.

Un $\varphi$-module sur $\calO$ est un $\calO$-module de type fini
$D$ muni d'un morphisme $\varphi$-semi-lin\'eaire $\varphi:D
\rightarrow D$. On \'ecrit $\varphi^*(D)$ pour le $\calO$-module
engendr\'e par $\varphi(D)$ dans $D$ et on dit que $D$ est \'etale si
$D=\varphi^*(D)$. 
Un $\varphi$-module sur $\calE$ est un $\calE$-espace vectoriel de
dimension finie $D$ muni d'un morphsime $\varphi$-semi-lin\'eaire
$\varphi:D \rightarrow D$. On dit que $D$ est \'etale si $D$ a un
$\calO$-r\'eseau stable par $\varphi$ et \'etale.
Un $(\varphi,\Gamma)$-module est un $\varphi$-module muni d'une
action de $\Gamma$ par des morphismes semi-lin\'eaires 
(par rapport \`a l'action de $\Gamma$ sur les coefficients) 
et commutant \`a $\varphi$.

Si $D$ est un $(\varphi,\Gamma)$-module \'etale sur $\calE$, et si
$\widehat{\calE^{\rm nr}}$ est l'anneau construit par Fontaine dans
\cite{F90}, alors $V(D) \= (\widehat{\calE^{\rm nr}} \otimes_{\calE}
D)^{\varphi=1}$ est un $L$-espace vectoriel de dimension
$\dim_{\calE}(D)$ muni d'une action lin\'eaire et continue de $\g$:
c'est une repr\'esentation $L$-lin\'eaire de $\g$. On a alors le
r\'esultat suivant (\cite{F90}) :

\begin{theo}\label{fontequiv}
Le foncteur $D \mapsto V(D)$ est une \'equivalence de cat\'egories 
\begin{itemize}
\item[(i)] de la cat\'egorie des $(\varphi,\Gamma)$-modules \'etales sur
$\calE$ vers la cat\'egorie des repr\'esentations $L$-lin\'eaires de
$\g$;
\item[(ii)] de la cat\'egorie des $(\varphi,\Gamma)$-modules \'etales sur
$\calO$ vers la cat\'egorie des $\O$-repr\'esenta\-tions de
$\g$.
\end{itemize}
\end{theo}

L'inverse de ce foncteur est not\'e $V \mapsto D(V)$. Les
$(\varphi,\Gamma)$-modules \'etales nous donnent donc un moyen
commode de travailler avec les repr\'esentations $L$-lin\'eaires. En
th\'eorie, les $(\varphi,\Gamma)$-modules permettent de retrouver
tous les objets que l'on peut associer aux repr\'esentations
$L$-lin\'eaires (et en pratique, c'est souvent le cas). Dans ce
paragraphe et les suivants, nous allons rappeler quelques unes de ces
constructions. 

Si $D$ est un $\varphi$-module \'etale et si $x \in D$, alors 
on peut \'ecrire $x=\sum_{i=0}^{p-1} (1+X)^i \varphi(x_i)$ o\`u les
$x_i \in D$ sont bien d\'etermin\'es.

\begin{defi}\label{defipsi}
On d\'efinit un op\'erateur $\psi : D \rightarrow
D$ par la formule $\psi(x) \= x_0$ si $x=\sum_{i=0}^{p-1} (1+X)^i
\varphi(x_i)$.
\end{defi}

Cet op\'erateur v\'erifie alors $\psi(\varphi(x)y) = x \psi(y)$ et
$\psi(x \varphi(y)) = \psi(x) y$ si $x \in \calO$ (ou $\calE$) 
et $y \in D$, et commute \`a l'action de $\Gamma$ si $D$ est un
$(\varphi,\Gamma)$-module.
C'est cet op\'erateur qui permet de faire le lien entre
les repr\'esentations $L$-lin\'eaires, les repr\'esentations de $\G$, et
la cohomologie d'Iwasawa, dont nous rappelons la
d\'efinition ci-dessous. Pour all\'eger les notations, on pose $F_n \=
\Q (\pmb{\mu}_{p^n})$ si $n \in \mathbb{N} \cup \{\infty\}$.

\begin{defi}\label{defiiw}
Si $V$ est une repr\'esentation $L$-lin\'eaire, 
si $T$ est un $\O$-r\'eseau $\g$-stable de $V$,
et $i \geq 0$, alors on d\'efinit :
$$ H^i_{\rm Iw}(\Q,T) \=
\varprojlim_{{\rm cor}_{F_{n+1}/F_n}} H^i({\rm
  Gal}(\Qpbar/F_n),T).$$  
\end{defi}
Le groupe $H^i_{\rm Iw}(\Q,V) \= 
\Q \otimes_{\Z} H^i_{\rm Iw}(\Q,T)$
ne d\'epend alors pas du choix d'un 
$\O$-r\'eseau $\g$-stable $T$ de $V$. Les 
$\Q \otimes_{\Z} \Z[[\Gamma]]$-modules $H^i_{\rm Iw}(\Q,V)$
ont \'et\'e \'etudi\'es en d\'etail par Perrin-Riou, qui a notamment
montr\'e (cf. \cite[\S 3.2]{BP94}) :

\begin{prop}\label{bpriw}
Si $V$ est une repr\'esentation $p$-adique de $\g$, alors on a $H^i_{\rm
Iw}(\Q,V)=0$ si $i \neq 1,2$. De plus : 
\begin{itemize}
\item[(i)] le sous-module de $\Q \otimes_{\Z} \Z[[\Gamma]]$-torsion de $H^1_{\rm
Iw}(\Q,V)$ est isomorphe \`a $V^{{\rm Gal}(\Qpbar / F_\infty)}$ 
et $H^1_{\rm Iw}(\Q,V)/\{{\rm torsion}\}$ est un $\Q \otimes_{\Z}
  \Z[[\Gamma]]$-module libre de rang $\dim_{\Q}(V)$;
\item[(ii)] $H^2_{\rm Iw}(\Q,V)=(V^*(1)^{{\rm Gal}(\Qpbar /
  F_\infty)})^*$. 
\end{itemize}
\end{prop}

Les espaces $D(V)^{\psi=1}$ et $D(V)/(1-\psi)$ sont eux aussi des $\Q
\otimes_{\Z} \Z[[\Gamma]]$-modules et on a le r\'esultat
suivant (cf. \cite[\S II.1]{CC99}) : 

\begin{prop}\label{cciw}
On a des isomorphismes canoniques de $\Q
\otimes_{\Z} \Z[[\Gamma]]$-modules
$H^1_{\rm Iw}(\Q,V) \simeq D(V)^{\psi=1}$ et $H^2_{\rm Iw}(\Q,V) \simeq
D(V)/(1-\psi)$. 
\end{prop}

En d'autres termes, la cohomologie du complexe $0 \rightarrow D(V)
\overset{1-\psi}{\longrightarrow} D(V) 
\rightarrow 0$ calcule la cohomologie d'Iwasawa de $V$.

\begin{coro}\label{dpsinozero}
Si $V$ est une repr\'esentation $L$-lin\'eaire, alors $D(V)^{\psi=1}
\neq 0$.
\end{coro}

Remarquons que le corollaire I.7.6 de \cite{CC99} affirme que
$D(V)^{\psi=1}$ contient en fait une base de $D(V)$ sur $\calE$. Nous
donnerons dans le corollaire \ref{iwnonzero2} une autre d\'emonstration du
corollaire \ref{dpsinozero}, qui ne passe pas par la cohomologie
d'Iwasawa (mais qui suppose que $V$ est cristalline). 

\subsection{Topologie faible et treillis}\label{topotr}

Nous allons maintenant nous int\'eresser 
\`a la topologie de $D$. Si $D$ est un
$\calO$-module libre de rang $d$, alors 
le choix d'une base de $D$ donne un isomorphisme
$D \simeq \calO^d$ et on peut
munir $D$ de la topologie faible. Un petit calcul montre qu'une
application $\calO$-lin\'eaire $\calO^d \rightarrow \calO^d$ est
n\'ecessairement continue et donc que la topologie d\'efinie sur $D$
par $D \simeq \calO^d$ ne d\'epend pas du choix d'une base de $D$. 

\begin{lemm}\label{bornetreillis}
Si $P$ est une partie d'un $\calO$-module libre $D$, 
et $M(P)$ est le $\O[[X]]$-module
engendr\'e par $P$, alors 
les propri\'et\'es suivantes sont \'equivalentes:
\begin{itemize}
\item[(i)] $P$ est born\'ee pour la topologie faible;
\item[(ii)] $M(P)$ est born\'e pour la topologie faible;
\item[(iii)] pour tout $j \geq 1$, l'image de $M(P)$ dans 
$D/\pi_L^j D$ est un $\O[[X]]$-module de type fini. 
\end{itemize}
\end{lemm}
\begin{proof}
Choisissons une base de $D$, et notons $D^+$ le $\O[[X]]$-module
engendr\'e par cette base, ce qui fait que la topologie faible sur $D$
est d\'efinie par le fait que les $\{\pi_L^i D + X^j D^+ \}_{i,j \geq 0}$ 
forment une base de voisinages de z\'ero. En particulier, une partie
$P$ de $D$ est born\'ee si et seulement si pour tout $k \geq 0$, il
existe $n(k,P) \in \mathbb{Z}$ 
tel que $P \subset \pi_L^k D + X^{n(k,P)} D^+$. Comme le
$\O[[X]]$-module engendr\'e par $\pi_L^k D + X^{n(k,P)} D^+$ est $\pi_L^k D +
X^{n(k,P)} D^+$ lui-m\^eme, on voit que les propri\'et\'es (i) et (ii)
sont \'equivalentes. Il reste donc \`a montrer que les $\O[[X]]$-modules born\'es sont
ceux qui satisfont (iii), c'est-\`a-dire 
qu'un $\O[[X]]$-module $M$ est born\'e si et seulement si pour
tout $j \geq 1$, l'image de $M$ dans $D/\pi_L^j D$ est un
$\O[[X]]$-module de type fini. Si $M$ est born\'e alors 
par d\'efinition, pour tout $j
\geq 1$, il existe $n(j,M)$ tel que $M \subset \pi_L^j D + X^{n(j,M)}
D^+$ ce qui fait que l'image de $M$ dans $D/\pi_L^j D$ est contenue
dans celle de $X^{n(j,M)} D^+$ et est de type fini puisque $\O[[X]]$
est un anneau noetherien. R\'eciproquement, si $M$ satisfait (iii), 
alors pour tout $j \geq 1$ l'image de $M$ dans $D/\pi_L^j D$ est un
$\O[[X]]$-module de type fini et est donc contenue dans $X^{n(j,M)}
D^+$ pour un $n(j,M) \in \mathbb{Z}$, ce qui fait que 
$M \subset \pi_L^j D + X^{n(j,M)} D^+$ et donc que $M$ est born\'e
pour la topologie faible.
\end{proof}

Un treillis de $D$ est un $\O[[X]]$-module born\'e $M \subset D$
tel que l'image de $M$ dans $D/\pi_L D$ en est un 
$k_L[[X]]$-r\'eseau. 
Un treillis d'un $\calE$-module $D$ est un treillis d'un
$\calO$-r\'eseau de $D$. 
Les treillis font l'objet d'une \'etude d\'etaill\'ee dans \cite[\S
4]{Co2}. Rappelons le r\'esultat de base sur les treillis 
stables par $\psi$ d'un $\varphi$-module $D$ (c'est la proposition
4.29 de \cite{Co2}):  

\begin{prop}\label{col429}
Si $D$ est un $\varphi$-module \'etale sur $\calO$, il existe un
unique treillis $D^{\sharp}$ de $D$ v\'erifiant les propri\'et\'es
suivantes :
\begin{itemize}
\item[(i)] quels que soient $x \in D$ et $k \in \N$, il existe $n(x,k) \in
  \N$ tel que $\psi^n(x) \in D^{\sharp} + p^k D$ si $n \geq n(x,k)$;
\item[(ii)] l'op\'erateur $\psi$ induit une surjection de $D^{\sharp}$ sur
  lui-m\^eme.
\end{itemize}
De plus :
\begin{itemize}
\item[(iii)] si $N$ est un sous-$\O[[X]]$ module born\'e 
de $D$ et $k \in \N$, il existe
$n(N,k)$ tel que $\psi^n(N) \subset D^{\sharp} + p^k D$ si $n \geq
n(N,k)$;
\item[(iv)] si $N$ est un treillis de $D$ stable par $\psi$ tel que
  $\psi$ induise une surjection de $N$ sur lui-m\^eme, alors $N
  \subset D^{\sharp}$ et $D^{\sharp} / N$ est annul\'e par $X$. 
\end{itemize}
\end{prop}

La fin de ce paragraphe est consacr\'e au d\'ebut de l'\'etude des
limites projectives d\'efinies ci-dessous.

\begin{defi}\label{defilimproj}
\begin{itemize}
\item[(i)] On note $(\varprojlim_{\psi} D(V))^{\rm b}$ le $L$-espace
  vectoriel des suites 
$(x_n)_{n \geq 0}$ d'\'el\'e\-ments de $D(V)$ telles que l'ensemble 
$\{x_n\}_{n \geq 0}$ est born\'e (d'o\`u le {\og b \fg}) pour 
la topologie faible et telles que $\psi(x_{n+1})=x_n$.
\item[(ii)] Si $T$ est un $\O$-r\'eseau de $V$ stable par $\g$, on note
$\varprojlim_{\psi}D^{\sharp}(T)$ le $\O$-module des suites
$\psi$-compatibles d'\'el\'ements de $D^{\sharp}(T)$, muni de 
la topologie de la limite projective.
\end{itemize}
\end{defi}

\begin{prop}\label{ddiese}
L'injection $D^{\sharp}(T) \hookrightarrow D(V)$
induit un isomorphisme topologique 
$L \otimes_{\O} \varprojlim_{\psi} D^{\sharp}(T)
\buildrel \sim \over \rightarrow (\varprojlim_{\psi}D(V))^{\rm b}$.  
\end{prop}
\begin{proof}
Il est clair que l'application ci-dessus est injective. 
Si $x=(x_n)_{n \geq 0} \in (\varprojlim_{\psi}D(V))^{\rm b}$, 
alors par d\'efinition
l'ensemble $P \= \{x_n\}_{n \geq 0}$ est born\'e 
pour la topologie faible
et par le lemme \ref{bornetreillis}, le $\O[[X]]$-module $M(P)$
engendr\'e par $P$ est born\'e (pour la topologie faible). 
Quitte \`a multiplier $x$ par une
puissance de $\pi_L$, on peut d'ailleurs supposer que $x_m \in D(T)$
pour tout $m \geq 0$. Si $k,m \geq 0$, alors
la proposition \ref{col429} appliqu\'ee \`a $N=M(P)$ montre que 
$x_m=\psi^n(x_{m+n}) \in D^{\sharp}(T) + p^k D(T)$ si $n$ est assez
grand. Comme c'est vrai pour tout $k$, on en d\'eduit que $x_m \in
D^{\sharp}(T)$ pour tout $m$ et donc que l'application 
$L \otimes_{\O} \varprojlim_{\psi}D^{\sharp}(T) \rightarrow   
(\varprojlim_{\psi}D(V))^{\rm b}$ est une bijection. C'est un
hom\'eomorphisme car la topologie de $D^{\sharp}(T)$ est la topologie
induite par la topologie faible de $D(T)$ via l'inclusion $D^{\sharp}(T)
\subset D(T)$.     
\end{proof}

Rappelons que $\psi:D^{\sharp}(T)\rightarrow D^{\sharp}(T)$ est
surjective, ce qui fait que les applications de transition dans 
$\varprojlim_{\psi}D^{\sharp}(T)$ le sont.
Le lemme ci-dessous et son corollaire
seront utilis\'es dans le paragraphe \ref{laouiw}. 

\begin{lemm}\label{psietddiese}
L'application naturelle $(\varprojlim_{\psi}D^{\sharp}(T)) / (1-\psi)
\rightarrow D^{\sharp}(T) / (1-\psi)$ est un isomorphisme.
\end{lemm}
\begin{proof}
Cette application est \'evidemment surjective, et nous allons montrer
qu'elle est injective, c'est-\`a-dire que si  $x=(x_n)_{n \geq 0} \in
\varprojlim_{\psi}D^{\sharp}(T)$, avec $x_0 \in
(1-\psi)D^{\sharp}(T)$, alors $x \in (1-\psi)
\varprojlim_{\psi}D^{\sharp}(T)$. Soit $y_0 \in D^{\sharp}(T)$ tel que
$(1-\psi)y_0 = x_0$. Pour tout $m \geq 0$, il existe $y^0_m \in
D^{\sharp}(T)$ tel que $\psi^m(y_m^0) = y_0$ et on a alors
$(1-\psi)y_m^0 - x_m \in D^{\sharp}(T)^{\psi^m =0}$. L'op\'erateur
$1-\psi$ est 
bijectif sur $D^{\sharp}(T)^{\psi^m =0}$ 
(un inverse \'etant donn\'e par $1+\psi+\psi^2+\cdots+\psi^{m-1}$)
et il existe
donc $y_m \in D^{\sharp}(T)$ tel que $(1-\psi)y_m = x_m$. Pour tout $k
\geq 0$, soit $z_k \= (z_{k,n})_{n \geq 0} \in \varprojlim_{\psi}
D^{\sharp}(T)$ tel que $z_{k,k} = y_k$.
Comme $\varprojlim_{\psi}D^{\sharp}(T)$ est
un espace topologique compact, la suite $\{z_k\}_{k \geq 0}$ a une
valeur d'adh\'erence $z$ et comme $(1-\psi)z_k \rightarrow x$ quand $k
\rightarrow \infty$ (par la continuit\'e de $\psi$, voir la
proposition \ref{actgcont} ci-dessous), 
on voit que $(1-\psi)z = x$ et donc que $x \in (1-\psi)
\varprojlim_{\psi}D^{\sharp}(T)$.
\end{proof}

\begin{coro}\label{h2ddiese}
On a un isomorphisme de $\O[[\Gamma]]$-modules
$(\varprojlim_{\psi}D^{\sharp}(T)) / (1-\psi) \simeq H^2_{\rm Iw}(\Q,T)$. 
\end{coro}
\begin{proof}
La proposition 4.43 de \cite{Co2}  
affirme que $D^{\sharp}(T) / (1-\psi) =
D(T)/(1-\psi)$, et  la remarque II.3.2 de \cite{CC99} 
nous dit  que
$D(T)/(1-\psi) = H^2_{\rm Iw}(\Q,T)$ ce qui fait que l'on a une
succession d'isomorphismes : 
$$(\varprojlim_{\psi}D^{\sharp}(T)) / (1-\psi) \simeq D^{\sharp}(T) /
(1-\psi) D^{\sharp}(T) \simeq  D(T)/(1-\psi)D(T) \simeq  H^2_{\rm
  Iw}(\Q,T). $$ 
\end{proof}

\section{Repr\'esentations cristallines de $\g$}\label{chapcris}

Le but de cette partie est de d\'ecrire explicitement $D^{\sharp}(V)$ en
fonction de $\dcris(V)$ lorsque $V$ est une repr\'esentation
$L$-lin\'eaire cristalline de dimension $2$ de $\g$ (et absolument 
irr\'eductible). Un \'enonc\'e similaire est d\'ej\`a paru dans \cite[\S4.4]{Co1}, mais sans preuve. Dans le dernier paragaphe de ce
chapitre, nous d\'efinissons une action d'un sous-groupe du Borel $\B$ sur $L \otimes_{\O} \varprojlim_{\psi} D^{\sharp}(T)$ et nous montrons
qu'elle est continue, topologiquement 
irr\'eductible, et que le lemme de Schur est
vrai pour cette action. 

\subsection{Rappels sur $\dcris(V)$}\label{cristallin}

Le but de ce paragraphe est de rappeler la classification des
repr\'esentations cristallines $\varphi$-semi-simples $V$ de dimension
$2$ de $\g$ en fonction de leur $\varphi$-module filtr\'e $\dcris(V)$. 

Afin de classifier certaines repr\'esentations $L$-lin\'eaires du groupe
$\g$, Fontaine a introduit les anneaux $\bcris$ et
$\bdR$. Ces anneaux v\'erifient les propri\'et\'es suivantes:
\begin{itemize}
\item[(i)] L'anneau $\bcris$ est une $\Q$-alg\`ebre munie d'une action
  de $\g$, telle que $\bcris^{\g}=\Q$ et d'un Frobenius $\varphi$ qui
  commute \`a l'action de $\g$;
\item[(ii)] le corps $\bdR$ est le corps des fractions d'un anneau 
  complet de valuation discr\`ete
  $\bdR^+$ (dont le corps r\'esiduel est
  $\C_p$), et il est donc muni de la filtration d\'efinie par les
  puissances de l'id\'eal maximal. Il est aussi muni d'une action
  continue de $\g$, telle que $\bdR^{\g}=\Q$ et la filtration 
  est donc stable sous l'action de $\g$;
\item[(iii)] on a une inclusion naturelle $\bcris \subset \bdR$ et une
  suite exacte (dite {\og fondamentale \fg}) : $$ 0 \rightarrow 
  \Q \rightarrow \bcris^{\varphi=1}
  \rightarrow \bdR/\bdR^+ \rightarrow 0;$$
\item[(iv)] il existe $t \in \bcris$ tel que $\varphi(t)=pt$ et $t$ est un
  g\'en\'erateur de l'id\'eal maximal de $\bdR^+$. 
  Le choix d'un tel $t$ d\'etermine une 
  application injective $\calR^+ \rightarrow L \otimes_{\Q} \bcris$ 
  telle que $t=\log(1+X)$, et cette injection commute \`a $\g$ et \`a
  $\varphi$; 
\item[(v)] si $m \geq 0$, alors on a une application toujours injective
  $\varphi^{-m} : L \otimes_{\Q}
  \bcris \rightarrow \bdR$ et la filtration de $\calR^+$ induite par
  cette application est donn\'ee par la filtration {\og ordre
  d'annulation en $\zeta_{p^m}-1$ \fg} o\`u $\zeta_{p^m}$ est une
  racine primitive $p^m$-i\`eme de $1$ d\'etermin\'ee par le choix de
  $t$.  
\end{itemize}

\'Etant donn\'ee une repr\'esentation $L$-lin\'eaire $V$, on
pose :
\begin{equation}\label{defdcrisddr}
\dcris(V) \= (\bcris \otimes_{\Q} V)^{\g}
\quad\text{et}\quad\ddR(V) \= (\bdR
\otimes_{\Q} V)^{\g}.
\end{equation} 

En g\'en\'eral, ces $\Q$-espaces vectoriels sont de
dimensions inf\'erieure ou \'egale \`a $\dim_{\Q}(V)$ et on dit que $V$ est
cristalline (resp. de de Rham) si $\dim_{\Q} \dcris(V)$
(resp. $\dim_{\Q} \ddR(V)$) est \'egale \`a $\dim_{\Q}(V)$. 
Comme $\bcris \subset \bdR$, une repr\'esentation cristalline
est n\'ecessairement de de Rham. Toutes les
repr\'esentations que nous consid\'erons dans cet article sont
suppos\'ees cristallines. 

Le Frobenius $\varphi$ de $\bcris$ commute \`a l'action de
$\g$ et la filtration de $\bdR$ est stable par $\g$, 
ce qui fait que $\dcris(V)$
est un $\varphi$-module et que $\ddR(V)$ est un $\Q$-espace vectoriel
filtr\'e. Si $V$ est une repr\'esentation cristalline, 
alors l'application naturelle de 
$\dcris(V)$ dans $\ddR(V)$ est un isomorphisme 
et $\dcris(V)$ est donc un $\varphi$-module filtr\'e. 
Si $V$ est $L$-lin\'eaire, alors $\dcris(V)$ et $\ddR(V)$ sont
naturellement des $L$-espaces vectoriels, et $\varphi:\dcris(V)
\rightarrow \dcris(V)$ ainsi que la filtration sur $\ddR(V)$ sont 
$L$-lin\'eaires.

Si $V$ est une repr\'esentation de de Rham, les poids de
Hodge-Tate de $V$ sont par d\'efinition les entiers $h$ tels que ${\rm Fil}^{-h}
\ddR(V) \neq {\rm Fil}^{-h+1} \ddR(V)$, compt\'es avec la
multiplicit\'e $\dim_L {\rm Fil}^{-h}
\ddR(V) / {\rm Fil}^{-h+1} \ddR(V)$.

Si $D$ est un $\varphi$-module de dimension $1$, on d\'efinit
$t_N(D) \= {\rm val}(\alpha)$ o\`u $\alpha$ 
est la matrice de $\varphi$ dans une
base de $D$ et si $D$ est un espace vectoriel filtr\'e de dimension
$1$, on d\'efinit $t_H(D)$ comme \'etant le plus grand $h \in
\mathbb{Z}$ tel 
que ${\rm Fil}^h D \neq 0$. Si $D$ est un $\varphi$-module de
dimension $\geq 1$, on d\'efinit $t_N(D) \= t_N(\det D)$ et $t_H(D) \=
t_H(\det D)$, ce qui fait que $t_H(D)$ est aussi l'oppos\'e de la
somme des poids de Hodge-Tate de $V$, compt\'es avec multiplicit\'es. 

Si $D$ est un $\varphi$-module filtr\'e, on dit que $D$ est
admissible si $t_N(D) = t_H(D)$ et si $t_N(D') - t_H(D') \geq
0$ pour tout sous-$\varphi$-module $D' \subset D$. 
Le fait que pour tout $h \geq 1$, on a ${\rm Fil}^{h+1}
\bcris^{\varphi=p^h} = 0$
permet de montrer que si $V$ est une repr\'esentation
cristalline, alors $\dcris(V)$ est admissible, et on a alors le
th\'eor\`eme de Colmez-Fontaine (voir \cite[th\'eor\`eme A]{CF} et
aussi \cite{Be2}) :

\begin{theo}\label{CFthm}
Le foncteur $\dcris(\cdot)$ est une \'equivalence de cat\'egories
de la cat\'ego\-rie des repr\'esentations $L$-lin\'eaires
cristallines dans la cat\'egorie des $\varphi$-modules filtr\'es
$L$-lin\'eaires admissibles.
\end{theo}

Ce th\'eor\`eme nous permet de faire une {\og liste \fg} des
repr\'esentations cristallines qui nous int\'eressent dans cet
article, qui sont les repr\'esentations cristallines $L$-lin\'eaires
$\varphi$-semi-simples
dont les poids de Hodge-Tate sont $0$ et $k-1$ avec $k \geq 2$. 
Soit $V$ une telle repr\'esentation. Quitte
\`a remplacer $L$ par un corps plus gros, on peut supposer que $L$
contient les valeurs propres $\alpha^{-1}$ et $\beta^{-1}$ de
l'op\'erateur $\varphi: \dcris(V) \rightarrow
\dcris(V)$. Quitte \`a \'echanger $\alpha$ et $\beta$, 
on peut supposer ${\rm val}(\alpha) \geq {\rm val}(\beta)$. 
L'hypoth\`ese de semi-simplicit\'e de cet op\'erateur 
implique que l'on peut \'ecrire
$\dcris(V)=L \cdot e_\alpha \oplus L \cdot e_\beta$ o\`u
$\varphi(e_\alpha) = \alpha^{-1} e_\alpha$ et $\varphi(e_\beta) =
\beta^{-1} e_\beta$. 
Il existe alors une droite $\Delta$ de $\dcris(V)$ telle que 
la filtration sur $\dcris(V)$ est donn\'ee par :
$${\rm Fil}^i \dcris(V) = \begin{cases}
\dcris(V) & \text{si $i \leq -(k-1)$;} \\
\Delta & \text{si $-(k-2) \leq i \leq 0$;} \\
0 & \text{si $i \geq 1$.}
\end{cases}$$

\`A isomorphisme pr\`es (et quitte \`a multiplier $e_\alpha$ ou  
$e_\beta$ par un \'el\'ement de $L$) cette droite est soit $L
\cdot e_\alpha$, soit $L \cdot e_\beta$, soit $L \cdot (e_\alpha +
e_\beta)$. Le fait que $\dcris(V)$ est admissible se traduit
par le fait que $0 \leq {\rm val}(\beta) \leq {\rm val}(\alpha) \leq
k-1$ et par ces conditions sur $\Delta$:
\begin{itemize}
\item[(i)] Si $0 < {\rm val}(\beta) \leq {\rm val}(\alpha) < k-1$, alors
$\alpha \neq \beta$, $\Delta = L \cdot (e_\alpha + e_\beta)$ et $V$
  est irr\'eductible;
\item[(ii)] si $0 = {\rm val}(\beta)$, alors $\Delta \neq L \cdot e_\alpha$ et :
\begin{itemize}
\item[1)] si $\Delta = L \cdot (e_\alpha +
  e_\beta)$, $V$ est r\'eductible non-scind\'ee;
\item[2)] si $\Delta = L \cdot e_\beta$, $V$ est r\'eductible scind\'ee.
\end{itemize}
\end{itemize} 

\begin{defi}\label{defdalbet}
Si $0 < {\rm val}(\beta) \leq {\rm val}(\alpha) < k-1$, on note
$D(\alpha,\beta)$ le $\varphi$-module filtr\'e admissible construit en (i)
ci-dessus.  
\end{defi}

\subsection{Modules de Wach des repr\'esentations cristallines}\label{rcm}
 
Dans ce paragraphe, nous faisons le lien entre la th\'eorie des
$(\varphi,\Gamma)$-modules et celle des repr\'esentations
cristallines de $\g$. L'id\'ee de d\'epart est que 
si $V$ est une repr\'esentation $p$-adique quelconque, alors son
$(\varphi,\Gamma)$-module $D(V)$ est un objet assez compliqu\'e,
mais dans le cas o\`u $V$ est cristalline, il est possible d'en
choisir une base particuli\`erement sympathique. Ce r\'esultat est
l'objet de la th\'eorie des modules de Wach (voir
\cite{F90,W96,Be1}). Rappelons le r\'esultat principal de la
th\'eorie (cf. \cite[\S II.1]{Be1}).

\begin{theo}\label{wachmod}
Si $V$ est une repr\'esentation $p$-adique de $\g$, alors $V$ est
cristalline \`a poids de Hodge-Tate dans un intervalle $[a;b]$
si et seulement s'il existe un sous-$L \otimes_{\O} \O[[X]]$-module
libre de rang fini $N(V) \subset D(V)$, tel que:
\begin{itemize}
\item[(i)] $D(V) = \calE \otimes_{L \otimes \O[[X]]} N(V)$;
\item[(ii)] $X^b N(V)$ est stable par $\varphi$ et $X^b N(V) / \varphi^*
  (X^b N(V))$ est annul\'e par $(\varphi(X)/X)^{b-a}$;
\item[(iii)] $N(V)$ est stable sous l'action de $\Gamma$ et celle-ci est
  triviale sur $N(V)/ X N(V)$. 
\end{itemize}
Si $T$ est un $\O$-r\'eseau de $V$, alors $N(T) \= D(T) \cap N(V)$ est un
$\O[[X]]$-module libre de rang fini qui satisfait l'analogue des
propri\'et\'es (i), (ii), (iii) ci-dessus.
\end{theo}

L'\'etude des modules de Wach fait l'objet d'une bonne partie de
\cite{Be1}. \`A partir de maintenant, on suppose que $V$ est une
repr\'esentation cristalline. L'objet de ce paragraphe est  
de d\'emontrer le th\'eor\`eme \ref{recupwach} ci-dessous. Avant de
l'\'enoncer, nous devons faire un certain nombre de rappels.

Il existe (cf. \cite[\S II.2]{Be1}) une 
application injective naturelle
(en particulier, $L \otimes_{\O} \O[[X]]$-lin\'eaire et commutant 
\`a $\varphi$, $\psi$ et \`a l'action de $\Gamma$) :
$$ N(V) \rightarrow \calR^+[1/t] \otimes_L \dcris(V), $$
application qui donne lieu \`a un isomorphisme (rappelons que
$t=\log(1+X)$) :
\begin{equation}\label{isomcriswach} 
\calR^+[1/t] \otimes_{L \otimes \O[[X]]} N(V) \buildrel\sim\over\rightarrow \calR^+[1/t]
\otimes_L \dcris(V). 
\end{equation}
Si les poids de Hodge-Tate de $V$ sont positifs, on a en fait : 
$$ N(V) \subset \calR^+ \otimes_L \dcris(V),$$
et la question se pose alors de caract\'eriser l'image de $N(V)$ dans $\calR^+ \otimes_L\dcris(V)$. C'est \`a cette question que r\'epond 
le th\'eor\`eme \ref{recupwach}.

Avant de donner cette caract\'erisation, rappelons que
pour tout $m \geq 0$, on a une application injective 
de $L$-alg\`ebres $\varphi^{-m} :
\calR^+ \rightarrow (\Q(\pmb{\mu}_{p^m}) \otimes_{\Q} L)[[t]]$, qui est
d\'etermin\'ee par le fait qu'elle envoie $\alpha \otimes X$ sur 
$\alpha \otimes (\zeta_{p^m} \exp(t/p^m) - 1)$ 
(ici $\zeta_{p^m}$ est une racine primitive
$p^m$-i\`eme de $1$ qui d\'epend du choix de $t\in\bcris$, voir le (v)
du d\'ebut du paragraphe \ref{cristallin}). 
Pour all\'eger les notations, on \'ecrit $L_m$ au lieu de
$\Q(\pmb{\mu}_{p^m}) \otimes_{\Q} L$. 
L'anneau $L_m[[t]]$ est muni d'une filtration par 
les puissances de l'id\'eal engendr\'e par
$t$, et cela d\'efinit une filtration sur $L_m[[t]] \otimes_L
\dcris(V)$.

Rappelons (cf. d\'efinition \ref{dordrer}) qu'un \'el\'ement $w\in \calR^+$
est d'ordre $r$ si, pour un (ou de mani\`ere \'equivalente tous les)
$0 < \rho <1$, la suite $(p^{-nr}\|w\|_{D(0,\rho^{1/p^n})})_n$ 
est born\'ee. On pose alors $\|w\|_r \= \sup_{n \geq 0}
p^{-nr}\|w\|_{D(0,\rho^{1/p^n})}$ (cette norme d\'epend bien s\^ur du choix
de $\rho$, mais pas la topologie qu'elle d\'efinit). 
On renvoie au lemme \ref{ordre} ci-dessous pour une 
caract\'erisation \'equivalente de la propri\'et\'e {\og \^etre
d'ordre $r$ \fg}. On a alors le r\'esultat suivant, qui
pr\'ecise le th\'eor\`eme D de \cite{Be2} (pour simplifier, on suppose
que $\varphi$ est semi-simple sur $\dcris(V)$) : 

\begin{theo}\label{recupwach}
Soit $V$ une repr\'esentation cristalline dont les poids de
Hodge-Tate sont positifs ou nuls et $\dcris(V)=\oplus_{i=1}^d L \cdot e_i$ avec
$\varphi(e_i)=\alpha_i^{-1} e_i$ et $\alpha_i \in L$. Si $x =
\sum_{i=1}^d x_i\otimes e_i \in \calR^+ 
\otimes_L \dcris(V)$, alors $x \in N(V)$ si et seulement si :
\begin{itemize}
\item[(i)] pour tout $1 \leq i \leq d$, $x_i$ est d'ordre $r_i \= 
{\rm val}(\alpha_i)$;
\item[(ii)] 
pour tout $m \geq 1$, $\sum_{i=1}^d \varphi^{-m}(x_i)\otimes \alpha_i^m
e_i \in {\rm Fil}^0 (L_m [[t]] \otimes_L \dcris(V))$.
\end{itemize} 
\end{theo}
\begin{proof}
Afin de d\'emontrer ce th\'eor\`eme, nous devons utiliser certains
r\'esultats et notations que nous n'avons pas rappel\'es dans cet article car 
ils ne sont pas utilis\'es ailleurs. Commen\c{c}ons par montrer que si $x \in N(V)$, alors il satisfait les
conditions (i) et (ii). On peut supposer que $x \in N(T)$ o\`u $T$ est
un r\'eseau $\g$-stable de $V$. Nous allons utiliser le fait que $N(T)$ est
naturellement un sous-ensemble de $\atplus[1/X] \otimes_{\Z} T$ o\`u
$\atplus$ est un sous-anneau de $\bcris^+$ tel que $\varphi:\atplus
\rightarrow \atplus$ est bijectif, et qui est born\'e pour la
topologie $p$-adique. En particulier, si $x \in N(T)$, alors
l'ensemble $\{\varphi^{-m}(x)\}_{m \geq 0}$ est born\'e et si 
$x = \sum_{i=1}^d x_i\otimes e_i$, alors $\varphi^{-m}(x)= 
\sum_{i=1}^d \varphi^{-m}(x_i)\otimes \alpha_i^m e_i$ ce qui fait que pour
tout $i \in \{1,\cdots,d\}$, la suite des $\{\varphi^{-m}(x_i)
\alpha_i^m\}_{m \geq 0}$ est born\'ee, et $x_i$ est donc d'ordre
$r_i = {\rm val}(\alpha_i)$ parce que si $\rho \geq p^{-1/(p-1)}$ et
$f(X) \in \calR^+$, alors $\|f(X)\|_{D(0,\rho)} =
\|\varphi(f(X))\|_{D(0,\rho^{1/p})}$. Cela montre le (i).
Le (ii) r\'esulte du fait que, comme $x \in \atplus[1/X] \otimes_{\Z} T$,
on a $\varphi^{-m}(x) \in \bdR^+ \otimes_{\Z} T$ 
pour tout $m \geq 1$ et donc : 
$$\varphi^{-m}(x) \in {\rm Fil}^0 (\bdR
\otimes_{\Z} T) \cap L_m [[t]] \otimes_L \dcris(V) = 
{\rm Fil}^0 (L_m [[t]] \otimes_L \dcris(V)).$$
Montrons maintenant que si $x \in \calR^+ \otimes_L \dcris(V)$
satisfait (i) et (ii), alors $x \in N(V)$. Pour cela, remarquons que
$x$ satisfait les conditions (1), (2) et (3) du th\'eor\`eme D de
\cite{Be2} (la condition (1) est automatique car $V$ est cristalline, la
condition (2) est \'equivalente \`a (ii) et la condition (3) \`a (i)). On en d\'eduit que, dans les notations de \cite{Be2}, $x \in
D^{\dagger}(V)$ o\`u $D^{\dagger}(V)$ est un certain sous-ensemble de
$D(V)$. Nous allons utiliser le fait que si $V$ est cristalline,
alors $D^{\dagger}(V) = \calE^{\dagger} \otimes_{L \otimes \O[[X]]}
N(V)$ o\`u $\calE^{\dagger}$ est le sous-corps de $\calE$
constitu\'e des s\'eries formelles $f(X) \in \calE$ ayant la
propri\'et\'e qu'il existe $0 \leq r(f) < 1$ tel que $f(X)$ 
converge sur la couronne $\{ z \in \C_p,\ r(f) \leq |z| < 1 \}$. Le fait que $x \in \calR^+ \otimes_L
\dcris(V)$ et que, comme on l'a rappel\'e plus haut, 
$\calR^+[1/t] \otimes_L \dcris(V) = 
\calR^+[1/t] \otimes_{L \otimes \O[[X]]}
N(V)$ montre que : $$x \in (\calR^+[1/t] \cap \calE^{\dagger})
\otimes_{L \otimes \O[[X]]} N(V).$$
Les coefficients de $x$ dans une base de $N(V)$ sont donc des
fonctions appartenant \`a $\calR^+[1/t] \cap \calE^{\dagger}$.
Si $f(X) \in \calR^+[1/\log(1+X)] \cap \calE^{\dagger}$, alors il
existe $\ell \geq 0$ tel que $f \in (L \otimes_{\O} \O[[X]])
[1/\varphi^{\ell}(X)]$. En effet, le fait que $f \in \calE^{\dagger}$ 
implique que $f$ ne peut avoir
qu'un nombre fini de p\^oles, et donc, puisque l'ensemble 
des z\'eros de $\log(1+X)$ est la r\'eunion de l'ensemble des z\'eros
des $\{ \varphi^{\ell}(X) \}_{\ell \geq 0}$, 
qu'il existe $\ell \geq 0$ tel que 
$f \in \calR^+ [1/\varphi^{\ell}(X)]$. Le fait que 
$f \in (L \otimes_{\O} \O[[X]]) [1/\varphi^{\ell}(X)]$
r\'esulte alors du fait que $\calR^+ \cap \calE^{\dagger}  = L
\otimes_{\O} \O[[X]]$. Il existe donc $\ell \geq 0$ tel
que $x \in N(V)[1/\varphi^\ell(X)]$. Le dernier r\'esultat que nous
allons utiliser sans d\'emonstration 
(mais qui suit des r\'esultats de \cite[\S 3.3]{Be1})
est que si $m \geq 1$, alors
$\varphi^{-m}(N(V))$ contient une base de ${\rm Fil}^0 (L_m[[t]]
\otimes_L \dcris(V))$ sur $L_m[[t]]$, ce qui fait que la condition
(ii) du th\'eor\`eme implique que les coefficients de $x$ dans une
base de $N(V)$ ne peuvent pas avoir de p\^oles ailleurs qu'en z\'ero,
et donc que $x \in N(V)[1/X]$. 
En effet, la filtration induite sur $\varphi^{-m}(\calR^+[1/t] \otimes
N(V))$ correspond \`a {\og l'ordre d'annulation \fg} en
$\zeta_{p^m}-1$, et le fait d'\^etre dans le ${\rm Fil}^0$ correspond
donc \`a l'absence de p\^ole. On termine en remarquant que les calculs de \cite[\S II.3 et \S III.2]{Be1} montrent que les diviseurs \'el\'ementaires de
l'inclusion : $$\calR^+ \otimes_{L \otimes \O[[X]]}
N(V) \subset \calR^+ \otimes_L \dcris(V)$$ sont premiers \`a $X$ et
donc que si $x \in \calR^+ \otimes_L \dcris(V)$, et $x \in N(V)[1/X]$, 
alors $x \in N(V)$.
\end{proof}

Soit $\prod_{i=1}^d \mathcal{H}^{r_i}$ l'espace de Banach des
$d$-uplets de fonctions holomorphes $(h_1,\cdots,h_d)$ o\`u $h_i \in
\calR^+$  est d'ordre $r_i$, muni de la norme $\|(h_1,\cdots, h_d)\|\=\sup_{1
\leq i \leq d} \|h_i\|_{r_i}$, et soit $N$ le sous-espace de
$\prod_{i=1}^d \mathcal{H}^{r_i}$ form\'e des fonctions qui satisfont
la condition (ii) du th\'eor\`eme \ref{recupwach}. Remarquons que $N(V)$ muni de la topologie $p$-adique est un espace de
Banach (on peut par exemple d\'ecider que $N(T)$ en est la boule
unit\'e). Nous aurons besoin du lemme suivant:

\begin{lemm}\label{istopon}
L'espace vectoriel $N$ est ferm\'e dans $\prod_{i=1}^d
\mathcal{H}^{r_i}$ et l'application naturelle 
$N(V) \rightarrow N$ est un
isomorphisme topologique d'espaces de Banach. 
\end{lemm}
\begin{proof}
Si $m \geq 1$, alors ${\rm Fil}^0 (L_m [[t]]
\otimes_L \dcris(V))$ est de codimension finie dans $L_m[[t]] 
\otimes_L \dcris(V)$ et le noyau de l'application :
$$\prod_{i=1}^d \mathcal{H}^{r_i} \rightarrow 
L_m [[t]] \otimes_L \dcris(V) / 
{\rm Fil}^0 (L_m [[t]] \otimes_L \dcris(V))$$ qui \`a
$(h_1,\cdots,h_d)$ associe 
l'image de $\sum_{i=1}^d 
\varphi^{-m}(h_i)\otimes\alpha_i^m e_i$ est donc un 
sous-espace ferm\'e $N_m$ de $\prod_{i=1}^d
\mathcal{H}^{r_i}$ ce qui fait que 
$N = \cap_{m \geq 1} N_m$ est lui
aussi ferm\'e dans $\prod_{i=1}^d
\mathcal{H}^{r_i}$. L'application naturelle $N(V) \rightarrow N$ est un
isomorphisme par le th\'eor\`eme
\ref{recupwach}.  Par le th\'eor\`eme de
l'image ouverte, pour montrer que $N(V) \rightarrow N$ est un
isomorphisme topologique, il suffit de montrer que l'application
$N(V) \rightarrow \prod_{i=1}^d \mathcal{H}^{r_i}$ est continue,
c'est-\`a-dire qu'elle est born\'ee. Pour cela, il suffit d'observer
qu'une boule unit\'e de $N(V)$ (pour la topologie $p$-adique) est
donn\'ee par $N(T)$, qui est un $\O[[X]]$-module de type fini.
\end{proof}

Le lemme \ref{istopon} montre que l'application $N(V) \rightarrow
\calR^+ \otimes_L \dcris(V)$ est continue pour la topologie forte. 

\begin{rema}\label{continj}
L'application $N(V) \rightarrow \calR^+ \otimes_L \dcris(V)$ est
continue si l'on donne \`a $N(V)$ la topologie faible et \`a
$\calR^+$ sa topologie de Fr\'echet.
\end{rema}
\begin{proof}
Cela r\'esulte du fait que l'application $L \otimes_{\O} \O[[X]]
\rightarrow \calR^+$ est continue, et du fait que $N(V)$ est un
$L \otimes_{\O} \O[[X]]$-module de type fini.
\end{proof}

\subsection{D'un monde \`a l'autre}\label{compare}

Le but de ce paragraphe est de donner une description explicite du
$L\otimes_{\O}\O[[X]]$-module $D^{\sharp}(V)$ en fonction du $L$-espace
vectoriel $\dcris(V)$. 
\'Etant donn\'es les calculs du paragraphe pr\'ec\'edent, il suffit de
faire le lien entre $N(V)$ et $D^{\sharp}(V)$.
On suppose tout au long de ce paragraphe que
$V$ est une repr\'esentation cristalline 
dont les poids de Hodge-Tate
sont dans un intervalle $[0;h]$. 
La proposition \ref{noslopezero} ci-dessous nous dit qu'alors $X^h N(V)
\subset D^{\sharp}(V) \subset N(V)$ (sous une petite condition
technique pour la deuxi\`eme inclusion). 

\begin{lemm}\label{wachstabpsi}
Si les poids de Hodge-Tate de $V$ sont positifs ou nuls et si $T$ est un $\O$-r\'eseau de $V$ stable par $\g$, alors 
les modules de Wach $N(V)$ et $N(T)$ sont stables par $\psi$. 
\end{lemm}
\begin{proof}
Rappelons que $h$ est un entier tel que les poids 
de Hodge-Tate de $V$ sont dans l'intervalle $[0;h]$.
Si $y \in N(T)$, alors par le (ii) du th\'eor\`eme \ref{wachmod}, 
on a $(\varphi(X)/X)^h (X^h y) \in \varphi^* (X^h N(T))$ et donc
$\varphi(X)^h y \in \varphi(X)^h \varphi^* (N(T))$ ce qui fait que $y
\in \varphi^* (N(T))$. Ceci \'etant vrai pour tout $y \in N(T)$, on a
donc $N(T) \subset \varphi^*(N(T))$ et donc $\psi(N(T)) \subset
N(T)$. Le fait que $\psi(N(V)) \subset N(V)$ en r\'esulte en inversant
$p$.  
\end{proof}

\begin{prop}\label{noslopezero}
Si les poids de Hodge-Tate de $V$ sont positifs ou nuls,
alors on a $X^h N(T) \subset D^{\sharp}(T) \subset X^{-1} N(T)$. Si les pentes de $\varphi$ sur $\dcris(V)$ sont de plus toutes strictement n\'egatives, 
alors on a $D^{\sharp}(T) \subset N(T)$.
\end{prop}
\begin{proof}
Comme $N(T)$ est un $\O[[X]]$-module qui engendre $D(T)$, 
c'est un treillis de $D(T)$ et 
le lemme \ref{wachstabpsi} montre que $N(T)$ est en plus
stable par $\psi$. Il suit alors de \cite[lemme 4.27, (iv) 
et proposition 4.29]{Co2} que $D^{\sharp}(T) \subset X^{-1} N(T)$. Par le th\'eor\`eme \ref{wachmod}, on a $\varphi(X^h N(T)) \subset X^h
N(T)$, ce qui fait que $X^h N(T) \subset \psi(X^h N(T))$. On sait par
ailleurs (par le lemme \ref{wachstabpsi}) que $\psi(N(T)) \subset
N(T)$, et on conclut de ces deux observations que la suite des
$\psi^n(X^h N(T))$ est une suite croissante de $\O[[X]]$-modules tous
contenus dans $N(T)$, ce qui fait qu'il existe $m(T) \geq 0$ tel que 
$\psi^{m+1}(X^h N(T)) = \psi^m(X^h N(T))$ pour tout $m \geq m(T)$, et
donc que $\psi^{m(T)}(X^h N(T))$ est un treillis qui contient
$X^h N(T)$ et sur lequel $\psi$ est surjectif. Le lemme r\'esulte alors du (iv) de la proposition \ref{col429}. Montrons finalement que si les pentes de $\varphi$ sur 
$\dcris(V)$ sont toutes $<0$, 
alors $D^{\sharp}(T) \subset N(T)$.
On vient de voir que $D^{\sharp}(T) \subset X^{-1} N(T)$. 
Pour terminer la d\'emonstration, nous allons 
montrer que pour tout $k \in \N$, 
il existe $n(k)$ tel que si $x \in N(T)$, alors $\psi^n(x) \in
N(T) + p^k X^{-1} N(T)$ si $n \geq n(k)$.
Le fait que $D^{\sharp}(T) = \psi^n(D^{\sharp}(T))$ implique alors que
$D^{\sharp}(T) \subset N(T) + p^k X^{-1} N(T)$, et comme c'est vrai pour
tout $k \geq 0$, c'est que $D^{\sharp}(T) \subset N(T)$. Nous allons donc montrer que $\psi^n \rightarrow 0$ sur le $L$-espace vectoriel 
de dimension finie $X^{-1} N(V) / N(V)$. 
Pour cela, rappelons que, par \cite[th\'eor\`eme
III.4.4]{Be1}, on a un isomorphisme de $\varphi$-modules $N(V)/X N(V)
= \dcris(V)$. Si $x \in X^{-1} N(V)$, alors on peut \'ecrire
$x=X^{-1} \sum_{i=1}^d \lambda_i \varphi(y_i)$ 
puisque $N(V) \subset \varphi^*
(N(V))$, ce qui fait que $\psi(x) = X^{-1} \sum_{i=1}^d \lambda_i(0) y_i
\mod{N(V)}$ et donc que, via l'identification $X^{-1} N(V)/ N(V)
= \dcris(V)$ qui \`a $\overline{X^{-1}x}$ 
associe $\overline{x}$, l'application
$\psi : X^{-1} N(V)/ N(V) \rightarrow X^{-1} N(V)/ N(V)$ correspond
\`a $\varphi^{-1} : \dcris(V) \rightarrow \dcris(V)$. Comme on a
suppos\'e que toutes les pentes de $\varphi$ sont $<0$, on a bien $\psi^n
\rightarrow 0$ sur $X^{-1} N(V)/ N(V)$,
ce qui fait que $D^{\sharp}(T) \subset N(T)$. 
\end{proof}

\begin{rema}\label{weakispt}
La topologie faible sur $D(T)$ induit une topologie sur $N(T)$ qui
n'est autre que la topologie $(p,X)$-adique, ce qui donne une
description plus agr\'eable de la topologie faible sur $D^{\sharp}(T)$. Ceci montre d'autre part qu'une partie de $N(V)$ est born\'ee pour la
topologie $p$-adique si et seulement si elle est born\'ee pour la
topologie faible. 
\end{rema}

\begin{coro}\label{iwnonzero2}
Les modules $D^{\sharp}(T)^{\psi=1}$ et donc $D(V)^{\psi=1}$ 
sont non nuls (comparer avec le corollaire \ref{dpsinozero}).  
\end{coro}
\begin{proof}
Il suffit de montrer que $N(T)^{\psi=1} \neq 0$. Si $x$ est non nul dans $N(T)$,
alors $y\=(1+X) \varphi(X^{h+1}x) \in (\varphi(X)N(T))^{\psi=0}$ comme
on le v\'erifie ais\'ement, et la s\'erie $\sum_{j=0}^{\infty}
\varphi^j(y)$ converge pour la topologie faible vers un \'el\'ement $z$ non nul dans $N(T)$ tel que $(1-\varphi)z=y$. On a donc $\psi(z)=z$ et $N(T)^{\psi=1} \neq 0$.
\end{proof}

Nous en venons maintenant au dernier point de ce paragraphe, la
description de $\varprojlim _{\psi} (D^{\sharp}(T))$ en termes de
fonctions \`a valeurs dans $\dcris(V)$.  
De l'application $N(V) \rightarrow \calR^+ \otimes_L \dcris(V)$, on
d\'eduit une application encore injective :
\begin{equation}\label{tjrsinj}
L \otimes_{\O} \varprojlim _{\psi} (D^{\sharp}(T)) 
\rightarrow \varprojlim _{\psi} \calR^+ \otimes_L 
\dcris(V).
\end{equation}
On suppose comme ci-dessus que $\varphi$ est semi-simple
sur $\dcris(V)$, et de plus
que $\varphi$ est de pentes strictement n\'egatives sur $\dcris(V)$. 

\begin{coro}\label{ddiesedcrisproj}
Soient $\dcris(V)=\oplus_{i=1}^d L \cdot e_i$ avec
$\varphi(e_i)=\alpha_i^{-1} e_i$ 
et $\alpha_i \in \pi_L\O$, et soient
$\{ (w_{i,n})_{n \geq 0} \}_{1\leq i\leq d}$ $d$ suites 
d'\'el\'ements de $\calR^+$ telles que 
$(v_n)_{n \geq 0} \= (\sum_{i=1}^d w_{i,n} \otimes e_i)_{n \geq 0} 
\in \varprojlim (\calR^+ \otimes_L \dcris(V))$. 
Alors $(v_n)_{n \geq 0} \in L \otimes_{\O} 
\varprojlim_{\psi}D^{\sharp}(T)$ si et 
seulement si : 
\begin{itemize}
\item[(i)] Pour tout $n \geq 0$ et tout $1 \leq i \leq d$, 
la fonction $w_{i,n}$ est d'ordre ${\rm val}(\alpha_i)$
et $\|w_{i,n}\|_{{\rm val}(\alpha_i)}$ est born\'e
ind\'ependamment de $n$;
\item[(ii)] pour tous $n\geq 0$, $m\geq 1$, on a :
$$\sum_{i=1}^d \varphi^{-m}(w_{i,n}\otimes e_i) \in {\rm
  Fil}^0(L_m[[t]] \otimes_L \dcris(V));$$  
\item[(iii)] pour tout $n\geq 1$ et tout $1 \leq i \leq d$,
$\psi(w_{i,n})=\alpha_i^{-1}w_{i,n-1}$.
\end{itemize} 
\end{coro}
 
\begin{proof} 
Par le th\'eor\`eme \ref{recupwach},
les conditions (i) et (ii) entra\^\i nent $v_n \in N(V)$. De plus, le fait que
$\|w_{i,n}\|_{{\rm val}(\alpha_i)}$ est born\'e 
ind\'ependamment de $n$ est alors \'equivalent, 
par le lemme \ref{istopon}, \`a ce que
la suite $(v_n)_{n \geq 0}$ est born\'ee dans $N(V)$ pour la topologie
$p$-adique, et donc aussi pour la topologie
faible (cf. remarque \ref{weakispt}). 
Enfin, la condition (iii) est \'equivalente
\`a ce que $\psi(v_{n+1})=v_n$ ce
qui \'equivaut \`a $(v_n)_n \in  (\varprojlim_{\psi} N(V))^{\rm b}$ et 
donc finalement \`a $(v_n)_n \in 
L \otimes_{\O} \varprojlim_{\psi} D^{\sharp}(T)$ par la proposition
\ref{ddiese}.  
\end{proof}

\subsection{Une action du Borel sup\'erieur}\label{borlim}

Le but de ce paragraphe est de d\'efinir une action du groupe : 
$$ G \= \left\{ \begin{pmatrix} 1 & \Q \\ 0 & \Q^{\times} 
\end{pmatrix} \right\} \subset \G $$ 
sur $\varprojlim_{\psi} D^{\sharp}(T)$ et de montrer
que cette action est continue et topologiquement irr\'eductible et 
que le lemme de Schur est v\'erifi\'e.
On suppose toujours que $V$ est une repr\'esentation cristalline
absolument irr\'eductible dont les poids de Hodge-Tate sont dans un
intervalle $[0;h]$. 
 
Tout \'el\'ement de $g \in G$ peut s'\'ecrire :
$$ g = \begin{pmatrix} 1 &  0 \\ 0 & p^k  
\end{pmatrix} \cdot \begin{pmatrix} 1 &  0 \\ 0 & a  
\end{pmatrix} \cdot \begin{pmatrix} 1 &  z \\ 0 & 1  
\end{pmatrix}, $$
o\`u $k\in\mathbb{Z}$, $a \in \Z^{\times}$ et $z \in \Q$. Si
$v=(v_n)_{n\geq 0} \in \varprojlim_{\psi}D^{\sharp}(T)$, on
\'etend la fonction $n \mapsto v_n$ \`a $\mathbb{Z}$ en posant $v_n \=
\psi^m(v_{n+m})$ pour $n+m\geq 0$. Si $a \in \Z^{\times}$, on note $[a] \in \Gamma$ l'unique \'el\'ement tel que $\varepsilon([a])=a$.  

\begin{defi}\label{actg}
Si $v=(v_n)_{n\geq 0} \in \varprojlim_{\psi}D^{\sharp}(T)$, on pose pour $n\in \N$ :
\begin{align*}
\left( \begin{pmatrix} 1 &  0 \\ 0 & p^k  
\end{pmatrix} \cdot v \right)_n & \= v_{n-k} = \psi^k(v_n) \\
\left( \begin{pmatrix} 1 &  0 \\ 0 & a  
\end{pmatrix} \cdot v \right)_n & \= [a](v_n) \\
\left( \begin{pmatrix} 1 &  z \\ 0 & 1  
\end{pmatrix} \cdot v \right)_n & \= \psi^m((1+X)^{p^{n+m} z} v_{n+m}),\ n+m\geq -{\rm val}(z).
\end{align*}
\end{defi}

On laisse le soin au lecteur de v\'erifier que les formules ci-dessus
d\'efinissent bien une action du groupe $G$ sur $\varprojlim_{\psi}D^{\sharp}(T)$. On d\'efinit aussi une structure de $\O[[X]]$-module sur
$\varprojlim_{\psi}D^{\sharp}(T)$ en posant $(1+X)^z \cdot v \=
\tiny{\begin{pmatrix} 1 & z \\ 0 & 1 \end{pmatrix}} \cdot v$ pour $z
\in \Z$. 

\begin{prop}\label{actgcont}
L'application $G \times \varprojlim_{\psi}D^{\sharp}(T) \rightarrow
\varprojlim_{\psi}D^{\sharp}(T)$ est continue. 
\end{prop}
\begin{proof}
On v\'erifie qu'il suffit de montrer que l'application $\psi:
\varprojlim_{\psi}D^{\sharp}(T) \rightarrow
\varprojlim_{\psi}D^{\sharp}(T)$ est continue et que l'application
$\BK \times \varprojlim_{\psi}D^{\sharp}(T)  \rightarrow
\varprojlim_{\psi}D^{\sharp}(T)$ est continue. Si $E$ et $\{X_i\}_{i\in I}$ sont des espaces topologiques
et si pour tout $i$, on se donne une application
continue $f_i : E \times X_i \rightarrow E \times X_i$, alors
l'application $E \times \prod_{i \in I} X_i \rightarrow \prod_{i \in
I} X_i$ donn\'ee par $(e,(x_i)_i) \mapsto (f_i(e,x_i))_i$ est
continue. En effet, la diagonale $\Delta_E$ de $\prod_{i \in I}
E$ y est ferm\'ee, et l'application  $(e,(x_i)_i) \mapsto (f_i(e,x_i))_i$
est la composition des applications :
$$E \times \prod_{i \in I} X_i = \Delta_E \times \prod_{i \in I} X_i
\subset \prod_{i \in I} (E \times X_i) \overset{\prod_{i \in I} f_i} 
{\longrightarrow}  \prod_{i \in I} (E \times X_i) 
\overset{\prod_{i \in I} {\rm pr}_i} {\longrightarrow}
\prod_{i \in I} X_i.  $$
On se ram\`ene donc \`a montrer que si chaque $f_i$
est continue, alors $\prod_{i \in I} f_i$ est continue (pour la
topologie produit) ce qui est laiss\'e en exercice facile au lecteur. Afin de montrer la proposition, il suffit donc de 
montrer que l'application $\psi:
D^{\sharp}(T) \rightarrow
D^{\sharp}(T)$ est continue, et que l'application
$\BK \times D^{\sharp}(T)  \rightarrow D^{\sharp}(T)$ est continue. Commen\c{c}ons par montrer que 
si $V$ est une repr\'esentation cristalline, alors 
(a) l'ensemble $\{p^j D(T) + X^k D^{\sharp}(T)\}_{j,k \geq 0}$ est
une base de voisinages de z\'ero pour la topologie faible et (b) l'ensemble $\{p^j D(T) + \varphi(X)^k D^{\sharp}(T)\}_{j,k \geq 0}$
est aussi une base de voisinages de z\'ero pour la topologie faible.
La proposition \ref{noslopezero} montre que 
$X^h N(T) \subset D^{\sharp}(T) \subset N(T)$,
ce qui montre le point (a) puisque $N(T)$ est un $\O[[X]]$-module 
libre qui engendre $D(T)$ sur $\calO$.
Pour montrer le point (b), 
et comme $p^j D(T) + \varphi(X)^k D^{\sharp}(T) \subset p^j
D(T) + X^k D^{\sharp}(T)$, 
il suffit de montrer (par exemple) que pour tout $k
\geq 0$, il existe $\ell$ tel que $p^j D(T) + X^\ell D^{\sharp}(T) \subset
p^j D(T) + \varphi(X)^k D^{\sharp}(T)$. Pour cela, remarquons que si $m
\geq j$ est tel que $p^m \geq k$, alors 
$p^j D(T) + X^{p^{m+1}} D^{\sharp}(T) \subset 
p^j D(T) + \varphi(X)^{p^m} D^{\sharp}(T)
\subset p^j D(T) + \varphi(X)^k D^{\sharp}(T)$ et on peut donc prendre
$\ell=p^{m+1}$. Le fait que l'op\'erateur $\psi: D^{\sharp}(T) \rightarrow 
D^{\sharp}(T)$ est continu pour la
topologie faible r\'esulte alors du fait que
$\psi\left(p^j D(T) + \varphi(X)^k D^{\sharp}(T)\right) 
=  p^j D(T) + X^k D^{\sharp}(T)$. Montrons maintenant que
l'application naturelle $\BK \times D^{\sharp}(T)
\rightarrow D^{\sharp}(T)$ est continue (pour la topologie
faible). Comme 
$\BK$ agit par multiplication par des \'el\'ements de $\O[[X]]$ ou
par action de $\Gamma$, et que la topologie faible de $D^{\sharp}(T)$
est d\'efinie par une base de voisinages de $0$ qui sont des
$\O[[X]]$-modules stables par $\Gamma$, 
on voit pour chaque $g \in \BK$
l'action de $g$ sur $D^{\sharp}(T)$ est continue. Il reste donc \`a montrer que si $W$ est un voisinage de z\'ero dans
$D^{\sharp}(T)$, alors il existe un sous-groupe normal ouvert $U$ de
$\BK$ tel que $u(x)-x \in W$ pour tous $(u,x) \in U \times
D^{\sharp}(T)$. Cela r\'esulte des faits suivants :
\begin{itemize}
\item[(i)] si $a \in \Z$ et $n \geq 0$, 
alors $(1+X)^{p^n a}-1 \in (p,X)^{n+1}$;
\item[(ii)] si $\gamma \in \Gamma$ et $n \geq 0$,  
alors $(\gamma^{(p-1)p^{n-1} }-1)N(T) \subset (p,X)^n N(T)$,
\end{itemize}
que nous laissons en
exercices au lecteur.
\end{proof}

La proposition suivante est le {\og lemme de Schur \fg} pour la
repr\'esentation de $G$ que l'on vient de d\'efinir.

\begin{prop}\label{schur}
Toute application continue
$G$-\'equivariante :
$$f: \varprojlim_{\psi}D^{\sharp}(T)
\longrightarrow \varprojlim_{\psi}D^{\sharp}(T)$$
est scalaire, i.e. est la multiplication par un \'el\'ement de $L$.
\end{prop}
\begin{proof}
Notons ${\rm pr}:
\varprojlim_{\psi}D^{\sharp}(T) \rightarrow D^{\sharp}(T)$ 
l'application $v \mapsto v_0$. Commen\c{c}ons par montrer que si
$v=(v_n)_{n \geq 0}$, alors ${\rm pr} \circ f(v)$ 
ne d\'epend que de $v_0 = {\rm pr}(v)$. Soit
$K_n$ l'ensemble des $v \in \varprojlim_{\psi}D^{\sharp}(T)$ dont les
$n$ premiers termes sont nuls, ce qui fait que pour $n \geq 1$, $K_n$
est un sous-$\O[[X]]$-module ferm\'e et stable par $\psi$ et $\Gamma$ 
de $\varprojlim_{\psi}D^{\sharp}(T)$ et que $\psi(K_n)=K_{n+1}$. On en
d\'eduit que ${\rm pr} \circ f (K_n)$ est un 
sous-$\O[[X]]$-module ferm\'e et stable par $\psi$ et $\Gamma$ 
de $D^{\sharp}(T)$. Par \cite[lemme 4.56]{Co2}, et comme $V$ est absolument
irr\'eductible, il existe une constante $\alpha(n) \in \N \cup \{
\infty \}$ telle que 
${\rm pr} \circ f (K_n)  = p^{\alpha(n)} D^{\sharp}(T)$. Le fait que
$\psi(K_n)=K_{n+1}$ et que $\psi(D^{\sharp}(T))=D^{\sharp}(T)$
implique $\alpha(n)=\alpha(n+1)$. Enfin, on voit que $K_n
\rightarrow 0$ (en ce sens que tout voisinage ouvert de $0$ dans
$\varprojlim_{\psi}D^{\sharp}(T)$ contient $K_n$ pour $n \gg 0$), ce
qui fait que $\alpha(n) \rightarrow + \infty$, et donc finalement que
$\alpha(n) = +\infty$ pour tout $n \geq 1$, ce qui veut dire que 
${\rm pr} \circ f (K_1) \subset K_1$ et donc que ${\rm pr} \circ f(v)$ 
ne d\'epend que de $v_0$. Pour tout $w \in D^{\sharp}(T)$, soit $\widetilde{w}$ un \'el\'ement de 
$\varprojlim_{\psi}D^{\sharp}(T)$ tel que $\widetilde{w}_0 = w$. Les
calculs pr\'ec\'edents montrent que l'application 
$h:  D^{\sharp}(T) \rightarrow D^{\sharp}(T)$ donn\'ee par
$h(w) \= {\rm pr} \circ f (\widetilde{w})$ est bien d\'efinie, et
qu'elle est $\O[[X]]$-lin\'eaire et commute \`a $\psi$ et \`a l'action
de $\Gamma$. Par \cite[proposition 4.7]{Co2}, elle s'\'etend en une
application de $(\varphi,\Gamma)$-modules $h: D(V) \rightarrow D(V)$ qui est scalaire car $V$ est irr\'eductible. Comme $f(v)_n = h(\psi^{-n} v)$ car $f$ commute \`a $\psi^{-n}$, on en d\'eduit que $f$ est aussi scalaire. 
\end{proof}

On a \'egalement :

\begin{prop}\label{actgirred}
L'action de $\B$ sur $\varprojlim_{\psi}D^{\sharp}(T)$ est
topologiquement irr\'e\-ducti\-ble. 
\end{prop}

\begin{proof}
C'est la r\'eunion de 
la d\'emonstration du corollaire 4.59
et du (iii) de la remarque 5.5 de \cite{Co2} (puisque $V$ 
et donc $V^*$ est absolument irr\'eductible). 
\end{proof}

\section{Repr\'esentations cristallines irr\'eductibles de $\G$}

Le but de cette partie est de d\'efinir des espaces de Banach $p$-adiques $\Pi(V)$ munis d'une action lin\'eaire continue de $\G$ et d'en commencer l'\'etude. Les repr\'esentations $\Pi(V)$ sont associ\'ees aux repr\'esentations cristallines irr\'eductibles $\varphi$-semi-simples $V$ de dimension $2$ de $\g$.

\subsection{Fonctions de classe ${\mathcal C}^r$ et distributions d'ordre $r$}\label{analyse}

Le but de ce paragraphe est de rappeler sans preuve les d\'efinitions et \'enonc\'es (classiques) d'analyse $p$-adique utilis\'es dans la suite.

Soit $L$ une extension finie de $\Q$. Si $f:\Z\rightarrow L$ est une fonction quelconque, on pose pour $n$ entier positif ou nul :
$$a_n(f)\=\sum_{i=0}^n(-1)^i\binom{n}{i}f(n-i).$$

Soit $r$ un nombre r\'eel positif ou nul.

\begin{defi}\label{classecr} 
Une fonction $f:\Z\rightarrow L$ est de classe ${\mathcal C}^r$ si
$n^r|a_n(f)|\rightarrow 0$ dans ${\mathbb R}_{\geq 0}$ quand $n\rightarrow
+\infty$. 
\end{defi} 

Rappelons que $f$ est continue si et seulement si $a_n(f)$ tend $p$-adiquement vers $0$ quand $n$ tend vers l'infini, de sorte que les fonctions de classe ${\mathcal C}^0$ au sens de la d\'efinition \ref{classecr} sont pr\'ecis\'ement les fonctions continues sur $\Z$. Toute fonction de classe ${\mathcal C}^r$ est aussi de classe ${\mathcal C}^s$ pour $0\leq s\leq r$ et est donc en particulier toujours continue. On note ${\mathcal C}^r(\Z,L)$ le $L$-espace vectoriel des fonctions $f:\Z\rightarrow L$ de classe ${\mathcal C}^r$.

Toute fonction continue $f:\Z\rightarrow L$ s'\'ecrit (d\'eveloppement de Mahler) :
\begin{equation}\label{mahler}
f(z)=\sum_{n=0}^{+\infty}a_n(f)\binom{z}{n} 
\end{equation} 
o\`u $z\in \Z$ et $\binom{z}{0}\=1$, $\binom{z}{n}=\frac{z(z-1)\cdots (z-n+1)}{n!}$ si $n\geq 1$. De plus, $\|f\|_0\={\rm sup}_{z\in\Z}|f(z)|$ co\"\i ncide avec ${\rm sup}_{n\geq 0}|a_n(f)|$. On v\'erifie facilement que l'espace ${\mathcal C}^r(\Z,L)$ est un espace de Banach pour la norme $\|f\|_r\={\rm sup}_{n\geq 0}(n+1)^r|a_n(f)|$.

La terminologie {\og de classe ${\mathcal C}^r$ \fg} provient du fait que, lorsque $r$ est entier strictement positif, les fonctions de ${\mathcal C}^r(\Z,L)$ sont aussi les fonctions sur $\Z$ qui, moralement, admettent $r$ d\'eriv\'ees avec la $r$-i\`eme d\'eriv\'ee continue, ce qui explique la terminologie. Par exemple, les fonctions localement analytiques sur $\Z$ sont de classe ${\mathcal C}^r$ pour tout $r\in \R_{\geq 0}$.

\begin{theo}[Amice-V\'elu, Vishik]\label{amice}
Soit $d$ un entier tel que $r-1<d$. Le sous-$L$-espace vectoriel ${\rm Pol}^d(\Z,L)$ de ${\mathcal C}^r(\Z,L)$ des fonctions $f:\Z\rightarrow L$ localement polynomiales de degr\'e (local) au plus $d$ est dense dans ${\mathcal C}^r(\Z,L)$.
\end{theo}

\begin{defi}
On appelle distribution temp\'er\'ee d'ordre $r$ un \'el\'ement de
l'espace ${\mathcal C}^r(\Z,L)^*$, c'est-\`a-dire une forme lin\'eaire continue sur l'espace de Banach des fonctions de classe ${\mathcal C}^r$.
\end{defi}

On dit parfois aussi distribution temp\'er\'ee d'ordre $\leq r$. Nous donnons maintenant deux descriptions des distributions temp\'er\'ees d'ordre $r$. 

Par le th\'eor\`eme \ref{amice}, l'inclusion ${\rm Pol}^d(\Z,L)\subsetneq {\mathcal C}^r(\Z,L)$ induit lorsque $r-1<d$ une injection :
$${\mathcal C}^r(\Z,L)^*\hookrightarrow {\rm Pol}^d(\Z,L)^*$$
o\`u ${\rm Pol}^d(\Z,L)^*$ est l'espace vectoriel dual de ${\rm Pol}^d(\Z,L)$. 
 
\begin{theo}[Amice-V\'elu, Vishik]\label{amice2}
Soit $\mu\in {\rm Pol}^d(\Z,L)^*$ et supposons que $r-1<d$. Alors $\mu\in {\mathcal C}^r(\Z,L)^*$ si et seulement s'il existe une constante $C_{\mu}\in L$ telle que, $\forall\ a\in \Z$, $\forall\ j\in \{0,\cdots,d\}$ et $\forall\ n\in {\mathbb N}$ :  
\begin{equation}\label{borneamice}
\int_{a+p^n\Z}(z-a)^jd\mu(z)\in C_{\mu}p^{n(j-r)}\O
\end{equation}
o\`u $\int_{a+p^n\Z}(z-a)^jd\mu(z)\=\mu\big({\bf 1}_{a+p^n\Z}(z)(z-a)^j\big)$ (${\bf 1}_{a+p^n\Z}$ est la fonction caract\'eristi\-que de $a+p^n\Z$).
\end{theo}

Notons que le plus petit entier $d$ tel que le th\'eor\`eme \ref{amice2} s'applique est la partie enti\`ere de $r$. Lorsque $\mu$ est d'ordre $r$ et $r-1<d$, on pose :
\begin{equation}\label{normedistrd}
\|\mu\|_{r,d}\={\rm sup}_{a\in\Z}{\rm sup}_{j\in \{0,\cdots,d\}}{\rm sup}_{n\in{\mathbb N}}\left(p^{n(j-r)}\left|\int_{a+p^n\Z}(z-a)^jd\mu(z)\right|\right).
\end{equation}
On peut montrer que $\|\mu\|_{r,d}$ est une norme sur ${\mathcal C}^r(\Z,L)^*$ qui redonne la topologie d'espace de Banach de ${\mathcal C}^r(\Z,L)^*$ et qui est \'equivalente \`a la norme :
\begin{equation}\label{normedistr}
\|\mu\|_{r}\={\rm sup}_{a\in\Z}{\rm sup}_{j,n\in{\mathbb N}}\left(p^{n(j-r)}\left|\int_{a+p^n\Z}(z-a)^jd\mu(z)\right|\right).
\end{equation}
En particulier, la majoration (\ref{borneamice}) est \'equivalente \`a la m\^eme majoration pour tout $j\in \N$ (et tout $a\in \Z$, $n\in \N$), quitte peut-\^etre \`a modifier $C_{\mu}$.

Nous aurons besoin du lemme suivant :

\begin{lemm}\label{lemmech}
Soit $r\in \R_{\geq 0}$ et $d$ la partie enti\`ere de $r$. Soit $n\in \N$ et :
$$f(z)\=\sum_{a\in \{0,\cdots,p^n-1\}}{\bf 1}_{a+p^n\Z}(z)\sum_{i=0}^d\lambda_{a,i}(z-a)^i\in {\rm Pol}^d(\Z,L)$$
o\`u $\lambda_{a,i}\in L$. Alors :
\begin{equation}\label{arghh}
\sup_{\mu\in {\mathcal C}^r(\Z,L)^*}\frac{\left|\int_{\Z}f(z)d\mu(z)\right|}{\|\mu\|_{r,d}}=\sup_{a\in \{0,\cdots,p^n-1\}}\sup_{i\in \{0,\cdots,d\}}\left|\lambda_{a,i}\right|p^{n(r-i)}.
\end{equation}
\end{lemm}
\begin{proof}
Par (\ref{normedistrd}), on voit que le r\'eel de gauche est plus petit que celui de droite. Pour $(a,i)\in \{0,\cdots,p^n-1\}\times \{0,\cdots,d\}$, il n'est pas difficile de construire une forme lin\'eaire $\mu_{a,i}\in {\rm Pol}^d(\Z,L)^*$ \`a support dans $a+p^n\Z$ satisfaisant (\ref{borneamice}) telle que $\int_{a+p^n\Z}(z-a)^jd\mu_{a,i}(z)=0$ si $j\ne i$ ($j\in \{0,\cdots, d\}$), $\int_{a+p^n\Z}(z-a)^id\mu_{a,i}(z)=p^{n(i-r)}$ et $\|\mu_{a,i}\|_{r,d}=1$ (les d\'etails sont laiss\'es en exercice au lecteur). En particulier, $\frac{\left|\int_{a+p^n\Z}f(z)d\mu_{a,i}(z)\right|}{\|\mu_{a,i}\|_{r,d}}=\left|\lambda_{a,i}\right|p^{n(r-i)}$ est inf\'erieur au r\'eel de gauche. Comme cela est vrai pour tout $(a,i)\in \{0,\cdots,p^n-1\}\times \{0,\cdots,d\}$, on en d\'eduit le r\'esultat.
\end{proof}

L'espace vectoriel des fonctions localement analytiques sur $\Z$,
${\rm An}(\Z,L)$ (muni de sa topologie d'espace de Fr\'echet, cf. \cite[\S16]{Sc}), est {\it a fortiori} dense dans ${\mathcal C}^r(\Z,L)$ et on dispose donc aussi d'une injection continue entre duaux continus :
$${\mathcal C}^r(\Z,L)^*\hookrightarrow {\rm An}(\Z,L)^*.$$
La transform\'ee d'Amice-Mahler :
\begin{equation}\label{Yvette} 
\mu\longmapsto \sum_{n=0}^{+\infty}\mu\left(\binom{z}{n}\right)X^n 
\end{equation} 
induit un isomorphisme topologique entre ${\rm An}(\Z,L)^*$ et
$\r^+$. Rappelons que (voir \S\ref{series}) : 
$$\r^+\=\left\{\sum_{n=0}^{+\infty}a_nX^n\mid a_n\in L,\ {\rm lim}_n|a_n|\rho^n=0\ \forall\ \rho\in [0,1[\right\}$$
et que cet espace est
muni de la topologie d'espace de Fr\'echet induite par la collection
des normes $\|\cdot \|_{D(0,\rho)}={\rm sup}_n(|a_n|\rho^n)$ pour
$0<\rho<1$. Rappelons \'egalement (d\'efinition \ref{dordrer}) qu'un \'el\'ement $w\in \r^+$ est d'ordre $r$ si, pour un (ou de mani\`ere \'equivalente tous les) $\rho\in ]0,1[$, la suite $\left(p^{-nr}\|w\|_{D(0,\rho^{1/p^n})}\right)_n$ est born\'ee dans $\R_{\geq 0}$.

\begin{lemm}\label{ordre}
Soit $w=\sum_{n=0}^{+\infty}a_nX^n\in \r^+$.
\begin{itemize}
\item[(i)] Un \'el\'ement $w=\sum_{n=0}^{+\infty}a_nX^n\in \r^+$ est
d'ordre $r$ si et seulement si $\{n^{-r}|a_n|\}_{n \geq 1}$ 
est born\'e (dans $\R_{\geq 0}$) lorsque $n$ varie.
\item[(ii)] Les normes ${\rm sup}_n\left(p^{-nr}\|w\|_{D(0,\rho^{1/p^n})}\right)$ et ${\rm sup}_n\left((n+1)^{-r}|a_n|\right)$ sont \'equivalentes pour tout $\rho\in ]0,1[$.
\end{itemize}
\end{lemm}
\begin{proof}
Voir \cite[\S1.1.4]{Co2}.
\end{proof}

Le r\'esultat suivant d\'ecoule imm\'ediatement des d\'efinitions et du lemme \ref{ordre} :  

\begin{prop}\label{amice3}
Soit $\mu\in {\rm An}(\Z,L)^*$. Alors $\mu\in {\mathcal C}^r(\Z,L)^*$ si et seulement l'\'el\'ement $\sum_{n=0}^{+\infty}\mu\left(\binom{z}{n}\right)X^n$ est d'ordre $r$ dans $\r^+$.
\end{prop}

On peut montrer que $\|\mu\|_{r}'\={\rm sup}_n((n+1)^{-r}\left|\mu\left(\binom{z}{n}\right)\right|)$ est une norme sur ${\mathcal C}^r(\Z,L)^*$ qui redonne la topologie d'espace de Banach de ${\mathcal C}^r(\Z,L)^*$.

\subsection{D\'efinition de $\Pi(V)$}\label{definition1}

Le but de ce paragraphe est de donner une premi\`ere d\'efinition de $\Pi(V)$ comme espace fonctionnel.

Soit $V$ une repr\'esentation cristalline $\varphi$-semi-simple absolument irr\'eductible comme au \S\ref{cristallin} (rappelons que $V$ \`a pour poids de Hodge-Tate $(0,k-1)$ avec n\'ecessairement $k\geq 2$ sinon $V$ n'est pas absolument irr\'eductible). Les valeurs propres $\alpha^{-1}$ et $\beta^{-1}$ de $\varphi$ sur $\dcris(V)=D(\alpha,\beta)$ sont alors telles que $\alpha\ne \beta$, ${\rm val}(\alpha)>0$, ${\rm val}(\beta)>0$, ${\rm val}(\alpha)+{\rm val}(\beta)=k-1$ et ${\rm val}(\alpha)\geq{\rm val}(\beta)$ (quitte \`a \'echanger $\alpha$ et $\beta$).

Soit $B(\alpha)$ l'espace de Banach suivant. Son $L$-espace vectoriel sous-jacent est form\'e des fonctions $f:\Q\rightarrow L$ v\'erifiant les deux conditions :
\begin{itemize}
\item[(i)] $f\!\mid_{\Z}$ est une fonction de classe ${\mathcal C}^{{\rm val}(\alpha)}$;
\item[(ii)] $(\frac{\alpha p}{\beta})^{{\rm val}(z)}z^{k-2}f(1/z)\!\!\mid_{\Z-\{0\}}$ se prolonge sur $\Z$ en une fonction de classe ${\mathcal C}^{{\rm val}(\alpha)}$.
\end{itemize}
Comme espace vectoriel, on a donc :
\begin{equation}\label{f1f2} 
B(\alpha)\simeq {\mathcal C}^{{\rm val}(\alpha)}(\Z,L)\oplus 
{\mathcal C}^{{\rm val}(\alpha)}(\Z,L),\ f\mapsto f_1\oplus f_2 
\end{equation} 
o\`u, pour $z\in \Z$, $f_1(z)\=f(pz)$ et $f_2(z)\=(\alpha
p\beta^{-1})^{{\rm val}(z)}z^{k-2}f(1/z)$. On munit $B(\alpha)$ de la norme :  
$$\|f\|\={\rm max}\big(\|f_1\|_{{\rm val}(\alpha)},\|f_2\|_{{\rm val}(\alpha)}\big),$$
qui en fait un espace de Banach en vertu du \S\ref{analyse}. On fait agir $L$-lin\'eairement $\G$ \`a gauche sur les fonctions de $B(\alpha)$ comme suit : 
\begin{equation}\label{action} 
\begin{pmatrix}a&b\\c&d\end{pmatrix}\!\cdot f(z)
=\alpha^{\!-{\rm val}(ad-bc)}\left(\frac{\alpha p}{\beta}\right)^{\!{\rm val}(-cz+a)}(-cz+a)^{\!k-2}
\!f\left(\frac{dz-b}{-cz+a}\right). 
\end{equation} 
Notons que {\tiny$\begin{pmatrix}a &0\\0&a\end{pmatrix}$} agit par la multiplication par $\varepsilon^{k-2}(a)(\frac{p^{k-1}}{\alpha\beta})^{{\rm val}(a)}\in \O^{\times}$.

\begin{lemm}\label{action2}
Si $f\in B(\alpha)$ et $g\in \G$, alors $g\cdot f\in B(\alpha)$ et l'action de $\G$ sur $B(\alpha)$ se fait par automorphismes continus.
\end{lemm}
\begin{proof}
On pose $r\= {\rm val}(\alpha)$, $d$ la partie enti\`ere de $r$ (donc $d\leq k-2$) et on munit l'espace de Banach ${\mathcal C}^r(\Z,L)$ de la norme induite par son bidual, c'est-\`a-dire :
$$\|f\|\=\sup_{\mu\in {\mathcal C}^r(\Z,L)^*}\frac{\left|\int_{\Z}f(z)d\mu(z)\right|}{\|\mu\|_{r,d}}$$
qui redonne la topologie d'espace de Banach de ${\mathcal C}^r(\Z,L)$ par \cite[lemme 9.9]{Sc} et les r\'esultats du \S\ref{analyse}. L'assertion du lemme est triviale si $g$ est scalaire. Par la d\'ecomposition de Bruhat et le cas scalaire, on est r\'eduit \`a montrer la stabilit\'e et la continuit\'e pour les matrices $g$ de la forme {\tiny$\begin{pmatrix}0 &p\\1&0\end{pmatrix}$}, {\tiny$\begin{pmatrix}1 &0\\0&\lambda\end{pmatrix}$} et {\tiny$\begin{pmatrix}1 &\lambda\\0&1\end{pmatrix}$} ($\lambda\in \Q^{\times}$). Le premier cas est \'evident puisqu'il envoie $f=(f_1,f_2)\in B(\alpha)$ sur $(f_2,f_1)\in B(\alpha)$ (cf. (\ref{f1f2})). Quitte \`a conjuguer par {\tiny$\begin{pmatrix}0 &p\\1&0\end{pmatrix}$} et \`a multiplier par un scalaire convenable, on peut prendre $\lambda\in \Z-\{0\}$ dans le deuxi\`eme et $g$ envoie alors $(f_1(z),f_2(z))\in B(\alpha)$ sur $(f_1(\lambda z),f_2(\lambda z))$ (\`a multiplication par des scalaires pr\`es). Il suffit donc de montrer que l'application $f(\cdot)\mapsto f(\lambda\cdot)$ est bien d\'efinie et continue de ${\mathcal C}^r(\Z,L)$ dans ${\mathcal C}^r(\Z,L)$. Par le th\'eor\`eme \ref{amice}, il suffit de montrer qu'il existe $c\in |L^{\times}|$ tel que, si $f\in {\rm Pol}^d(\Z,L)\subset {\mathcal C}^r(\Z,L)$, alors $\|f(\lambda\cdot)\|\leq c\|f(\cdot)\|$. En \'ecrivant, pour un $n\in \N$ assez petit et des $\lambda_{n,a,i}\in L$ convenables :
$$f(z)=\sum_{a\in \{0,\cdots,p^n-1\}}{\bf 1}_{a+p^n\Z}(z)\sum_{i=0}^d\lambda_{a,i}(z-a)^i,$$
un calcul donne :
$$f(\lambda z)=\sum_{\substack{a\in \{0,\cdots,p^n-1\} \\
{\rm val}(a)\geq {\rm val}(\lambda)}}
{\bf 1}_{\frac{a}{\lambda}+p^{n-{\rm val}(\lambda)}\Z}(z)\sum_{i=0}^d\lambda_{a,i}\lambda^i\left(z-\frac{a}{\lambda}\right)^i.$$
On en d\'eduit gr\^ace \`a (\ref{normedistrd}) :
\begin{eqnarray*}
\|f(\lambda\cdot)\|&\leq &\sup_{\substack{a\in \{ 0,\cdots,p^n-1 \} \\ {\rm val}(a)\geq {\rm val}(\lambda)}}\sup_{i\in \{0,\cdots,d\}}\left|\lambda_{a,i}\lambda^i\right| p^{(n-{\rm val}(\lambda))(r-i)}\\
&\leq &\sup_{a\in \{0,\cdots,p^n-1\}}\sup_{i\in \{0,\cdots,d\}}\left|\lambda_{a,i}\right| p^{n(r-i)}p^{-r{\rm val}(\lambda)}\\
&\leq & \|f(\cdot)\|
\end{eqnarray*}
o\`u la derni\`ere in\'egalit\'e r\'esulte du lemme \ref{lemmech} et du fait que ${\rm val}(\lambda)\geq 0$. Passons au dernier cas. Quitte \`a conjuguer $g=$ {\tiny$\begin{pmatrix}1 &\lambda\\0&1\end{pmatrix}$} par un \'el\'ement convenable de {\tiny$\begin{pmatrix}1 &0\\0&\Q^{\times}\end{pmatrix}$}, on peut supposer $\lambda=p$ et $g$ envoie alors $(f_1,f_2)$ sur $(f_1(z+1),(1+pz)^{k-2}f_2(z/(1+pz)))$ (\`a multiplication par des scalaires pr\`es). Il suffit donc de montrer que les applications $f(\cdot)\mapsto f(\cdot+1)$ et $f(\cdot)\mapsto (1+pz)^{k-2}f(\cdot /(1+p\cdot ))$ sont bien d\'efinies et continues de ${\mathcal C}^r(\Z,L)$ dans ${\mathcal C}^r(\Z,L)$. Cela se v\'erifie par un argument analogue au pr\'ec\'edent avec un calcul utilisant (\ref{normedistrd}) et (\ref{normedistr}) dont on laisse les d\'etails au lecteur.
\end{proof}

La repr\'esentation $B(\alpha)$ doit \^etre pens\'ee comme une induite parabolique. Sans donner un sens formel \`a ce qui suit, on a un isomorphisme $\G$-\'equivariant :
$$\big({\rm Ind}_{\B}^{\G}{\rm nr}(\alpha^{-1})\otimes d^{k-2}{\rm nr}(p\beta^{-1})\big)^{{\mathcal C}^{{\rm val}(\alpha)}}\simeq B(\alpha)$$
o\`u l'espace de gauche est celui des fonctions $F:\G\rightarrow L$ qui sont de classe ${\mathcal C}^{{\rm val}(\alpha)}$ (oublions que nous n'avons pas d\'efini de telles fonctions dans ce cadre !) telles que :
\begin{equation}\label{equfonc}
F\left( \begin{pmatrix}a&b\\0&d\end{pmatrix}g\right) = \alpha^{-{\rm val}(a)}p^{{\rm val}(d)}\beta^{-{\rm val}(d)}d^{k-2}F(g)
\end{equation}
avec action de $\G$ donn\'ee par $(g\cdot F)(g')\=F(g'g)$. On passe de $F$ \`a une fonction $f\in B(\alpha)$ en posant :
\begin{equation}\label{fonction} 
f(z)\=F\left(\begin{pmatrix}0&1\\-1&z\end{pmatrix}\right).
\end{equation}
On d\'efinit de m\^eme $B(\beta)$ muni d'une action de $\G$ par automorphismes en \'echangeant partout $\alpha$ et $\beta$. 

Voici des exemples importants de fonctions dans $B(\alpha)$ :

\begin{lemm}\label{appartient}
Pour $0\leq j<{\rm val}(\alpha)$ et $a\in \Q$, les fonctions $z\mapsto z^j$ et les fonctions :
$$z\longmapsto \left(\frac{\alpha p}{\beta}\right)^{{\rm val}(z-a)}\!\!(z-a)^{k-2-j}$$
(prolong\'ees par $0$ en $a$) sont dans $B(\alpha)$. 
\end{lemm} 
\begin{proof} 
En faisant agir {\tiny$\begin{pmatrix}0 &1\\1&0\end{pmatrix}$} sur $z^j$, il suffit de traiter les deuxi\`emes fonctions. En faisant agir {\tiny$\begin{pmatrix}1 &a\\0&1\end{pmatrix}$}, il suffit m\^eme par le lemme \ref{action2} de traiter le cas $a=0$ et comme $z\mapsto z^j$ est clairement de classe ${\mathcal C}^{{\rm val}(\alpha)}$ sur $\Z$, il suffit de montrer que $z\mapsto f(z)\=(\alpha p\beta^{-1})^{{\rm val}(z)}z^{k-2-j}$ est de classe ${\mathcal C}^{{\rm val}(\alpha)}$ sur $\Z$. Soit $f_0$ la fonction nulle sur $\Z$ et, pour $n\in {\mathbb Z}$, $n>0$, posons $f_n(z)\=(\alpha
p\beta^{-1})^{{\rm val}(z)}z^{k-2-j}$ si ${\rm val}(z)<n$, $f_n(z)\=0$ sinon. La fonction $f_n$ est de classe ${\mathcal C}^{{\rm val}(\alpha)}$ sur $\Z$ puisqu'elle est localement polynomiale. Il suffit de montrer que $f_{n+1}-f_n\rightarrow 0$ dans ${\mathcal C}^{{\rm val}(\alpha)}(\Z,L)$ quand $n\rightarrow +\infty$ (car $\sum_{n=0}^{\infty}(f_{n+1}-f_n)=f\in {\mathcal C}^{{\rm val}(\alpha)}(\Z,L)$ puisque ${\mathcal C}^{{\rm val}(\alpha)}(\Z,L)$ est complet). Par \cite[lemme 9.9]{Sc}, il suffit de v\'erifier que $f_{n+1}-f_n\rightarrow 0$ dans $(B(\alpha)^*)^*$, i.e. que :  
$$\sup_{\mu\in {\mathcal C}^{{\rm val}(\alpha)}(\Z,L)^*}\frac{\left|\int_{\Z}(f_{n+1}(z)-f_n(z))d\mu(z)\right|}{\|\mu\|_{{\rm val}(\alpha)}}\longrightarrow 0\ {\rm quand}\ n\rightarrow +\infty.$$  
Mais : 
$$\int_{\Z}(f_{n+1}(z)-f_n(z))d\mu(z)\!=\!\!\left(\frac{\alpha p}{\beta}\right)^{\!n}
\!\left(\int_{p^n\Z}\!z^{k-2-j}d\mu(z)-\int_{p^{n+1}\Z}\!z^{k-2-j}d\mu(z)\right)$$
d'o\`u :
$$\left|\int_{\Z}(f_{n+1}(z)-f_n(z))d\mu(z)\right|\leq c{\|\mu\|_{{\rm val}(\alpha)}}p^{-n(2{\rm val}(\alpha)-k+2)}p^{-n(k-2-j-{\rm val}(\alpha))}$$
o\`u $c\={\max}(1,p^{-(k-2-j-{\rm val}(\alpha))})$ en utilisant (\ref{normedistr}) et en se rappelant que ${\rm val}(\alpha)+{\rm val}(\beta)=k-1$. Cela implique $|\int_{\Z}(f_{n+1}(z)-f_n(z))d\mu(z)|\leq c{\|\mu\|_{{\rm val}(\alpha)}}p^{n(j-{\rm val}(\alpha))}$ d'o\`u le r\'esultat puisque $j<{\rm val}(\alpha)$.  
\end{proof} 

On a un lemme analogue en \'echangeant $\alpha$ et $\beta$. On note $L(\alpha)$ l'adh\'erence dans $B(\alpha)$ du sous-$L$-espace vectoriel engendr\'e par les fonctions $z^j$ et $(\alpha p\beta^{-1})^{{\rm val}(z-a)}(z-a)^{k-2-j}$ pour $a\in \Q$ et $j$ entier, $0\leq j<{\rm val}(\alpha)$. De m\^eme, on note $L(\beta)$ l'adh\'erence dans $B(\beta)$ du sous-$L$-espace vectoriel engendr\'e par les fonctions $z^j$ et $(\beta p\alpha^{-1})^{{\rm val}(z-a)}(z-a)^{k-2-j}$ pour $a\in \Q$ et $j$ entier, $0\leq j<{\rm val}(\beta)$. Notons que, lorsque $\alpha = p\beta$, $L(\beta)$ est de dimension finie et s'identifie aux polyn\^omes en $z$ de degr\'e au plus $k-2$.

\begin{lemm}
Le sous-espace $L(\alpha)$ (resp. $L(\beta)$) est stable par $\G$ dans $B(\alpha)$ (resp. $B(\beta)$.
\end{lemm}
\begin{proof} 
Exercice.
\end{proof} 

\begin{defi}
On pose $\Pi(V)\=B(\alpha)/L(\alpha)$.
\end{defi}

Il s'agit encore d'un $L$-espace de Banach (avec la topologie quotient) muni d'une action de $\G$ par automorphismes continus. Nous allons voir que l'application $\G\times \Pi(V)\rightarrow \Pi(V)$ est continue, que $\Pi(V)$ est unitaire et que l'on a un morphisme continu $\G$-\'equivariant $\widehat I:B(\beta)/L(\beta)\rightarrow \Pi(V)$ qui est un isomorphisme lorsque $\alpha\ne p\beta$.

\subsection{Autre description de $\Pi(V)$}\label{definition2}

Le but de ce paragraphe est de donner une description plus intrins\`eque de $\Pi(V)=B(\alpha)/L(\alpha)$ pour en d\'eduire certaines propri\'et\'es de l'action de $\G$ (continuit\'e, unitarit\'e, entrelacements) peu \'evidentes sur la d\'efinition pr\'ec\'edente. On conserve les notations du \S\ref{definition1}.

Soit :
$$\pi(\alpha)\={\rm Sym}^{k-2}L^2\otimes_L{\rm Ind}_{\B}^{\G}{\rm nr}(\alpha^{-1})\otimes {\rm nr}(p\beta^{-1})$$ 
la repr\'esentation de $\G$ produit tensoriel de la repr\'esentation
alg\'ebrique ${\rm Sym}^{k-2}L^2$ par l'induite parabolique lisse
${\rm Ind}_{\B}^{\G}{\rm nr}(\alpha^{-1})\otimes {\rm
  nr}(p\beta^{-1})$ (c'est-\`a-dire l'espace des fonctions localement
constantes $h:\G\rightarrow L$ v\'erifiant une \'egalit\'e analogue
\`a (\ref{equfonc}) avec action \`a gauche de $\G$ par translation \`a
droite). On munit $\pi(\alpha)$ de l'unique topologie localement
convexe (au sens de \cite{Sc}) telle que les ouverts sont les
sous-$\O$-modules g\'en\'erateurs (sur $L$). La repr\'esentation
$\pi(\alpha)$ est dite localement alg\'ebrique (cf. l'appendice de
\cite{ST2}) et n'est autre que la repr\'esentation ${\rm Alg}(V)
\otimes_L {\rm Lisse}(V)$ de l'introduction.

On identifie ${\rm Sym}^{k-2}L^2$ \`a l'espace vectoriel des polyn\^omes $P(z)$ de degr\'e au plus $k-2$ \`a coefficients dans $L$ munis de l'action \`a gauche de $\G$ :
\begin{equation}\label{algebraic}
\begin{pmatrix}a&b\\c&d\end{pmatrix}\!\cdot P(z)
=(-cz+a)^{\!k-2}P\left(\frac{dz-b}{-cz+a}\right).
\end{equation}
Comme en (\ref{fonction}), on identifie ${\rm Ind}_{\B}^{\G}{\rm nr}(\alpha^{-1})\otimes {\rm nr}(p\beta^{-1})$ \`a l'espace vectoriel des fonctions $f:\Q\rightarrow L$ localement constantes telles que $(\frac{\alpha p}{\beta})^{{\rm val}(z)}f(1/z)$ se prolonge sur $\Q$ en une fonction localement constante avec action \`a gauche de $\G$ comme en (\ref{action}) mais sans le facteur $(-cz+a)^{k-2}$. On en d\'eduit une injection $\G$-\'equivariante continue :
\begin{equation}\label{commefonction}
\pi(\alpha)\hookrightarrow B(\alpha),\ P(z)\otimes f(z)\mapsto P(z)f(z).
\end{equation} 
Par le th\'eor\`eme \ref{amice}, l'image de $\pi(\alpha)$ est dense dans $B(\alpha)$. En particulier, on a une injection continue $B(\alpha)^*\hookrightarrow \pi(\alpha)^*$. On d\'efinit de m\^eme $\pi(\beta)$ et une injection $\G$-\'equivariante continue d'image dense $\pi(\beta)\hookrightarrow B(\beta)$. Ces injections induisent des applications $\G$-\'equivariantes continues $\pi(\alpha)\rightarrow B(\alpha)/L(\alpha)$ et $\pi(\beta)\rightarrow B(\beta)/L(\beta)$.

Si $\pi^0$ est un sous-$\O$-module g\'en\'erateur d'un $L$-espace vectoriel $\pi$, rappelons qu'on appelle compl\'et\'e de $\pi$ par rapport \`a $\pi^0$ l'espace de Banach :
$$\Pi\= (\varprojlim_n \pi^0/\pi_L^n\pi^0)\otimes_{\O}L.$$
On a un morphisme canonique d'image dense $\pi\rightarrow \Pi$ qui n'est pas injectif en g\'en\'eral (si $\pi^0=\pi$, on a $\Pi=0$). Le dual continu $\Pi^*$ de $\Pi$ s'identifie en tant qu'espace de Banach \`a ${\rm Hom}_{\O}(\pi^0,\O)\otimes_{\O}L$ (avec ${\rm Hom}_{\O}(\pi^0,\O)$ comme boule unit\'e). Si $\pi$ est un espace localement convexe tonnel\'e (cf. \cite[\S6]{Sc}) muni d'une action continue d'un groupe topologique localement compact $G$ telle que $\pi^0$ est ouvert et stable par $G$, il est facile de v\'erifier en utilisant le th\'eor\`eme de Banach-Steinhaus (cf. \cite[proposition 6.15]{Sc}) que $\Pi$ et $\Pi^*$ sont des $G$-Banach unitaires et que la fl\`eche canonique $\pi\rightarrow \Pi$ est continue.

\begin{theo}\label{complete}
L'application $\pi(\alpha)\rightarrow B(\alpha)/L(\alpha)$ induit un isomorphisme topologique $\G$-\'equivariant entre $B(\alpha)/L(\alpha)$ et le compl\'et\'e de $\pi(\alpha)$ par rapport \`a un quelconque sous-$\O$-module g\'en\'erateur de $\pi(\alpha)$ stable par $\G$ et de type fini comme $\O[\G]$-module. On a le m\^eme r\'esultat en rempla\c cant $\alpha$ par $\beta$.
\end{theo}
\begin{proof}
Notons que le compl\'et\'e ne d\'epend pas du choix du sous-$\O[\G]$-module de type fini g\'en\'erateur de $\pi(\alpha)$ car ces $\O$-modules sont tous commensurables dans $\pi(\alpha)$. En utilisant $\G=\B\K$ et le fait que $\K$ est compact, on voit facilement qu'il suffit de compl\'eter par rapport \`a un sous-$\O[\B]$-module de type fini g\'en\'erateur quelconque, par exemple :
$$\sum_{j=0}^{k-2}\O[\B]({\mathbf 1}_{\Z}(z)z^j)+\sum_{j=0}^{k-2}
\O[\B]\left(\left(\frac{p\alpha}{\beta}\right)^{\!\!{\rm val}(z)}\!\!\!\!\!\!\!\!{\mathbf 1}_{\Q-\Z}(z)z^j\right)\subset \pi(\alpha)$$ 
o\`u ${\bf 1}_U$ est la fonction caract\'eristique de l'ouvert $U$. Le dual du compl\'et\'e cherch\'e est donc isomorphe au Banach :
\begin{multline}\label{thedual} 
\{\mu\in \pi(\alpha)^*\mid\forall\ g\in\B, \forall\ j\in
\{0,\cdots,k-2\}, |\mu(g({\mathbf 1}_{\Z}(z)z^j))|\leq 1\\   
 {\rm et}\ \left|\mu(g({\mathbf 1}_{\Q-\Z}(z)(\alpha p\beta^{-1})^{{\rm val}(z)}\!z^j))\right|\leq 1\}\otimes_{\O}L.  
\end{multline}
En utilisant l'int\'egralit\'e du caract\`ere central, il est \'equivalent de prendre $g\in $ {\tiny$\begin{pmatrix}1&\Q \\0&\Q^{\times}\end{pmatrix}$} dans
(\ref{thedual}). Pour $f\in \pi(\alpha)$, vue comme fonction sur $\Q$ via (\ref{commefonction}), et $U$ ouvert de $\Q$, on \'ecrit $\int_{U}f(z)d\mu(z)$ pour $\mu({\bf 1}_U(z)f(z))$. Un calcul donne alors que les conditions sur $\mu$ dans (\ref{thedual}) sont \'equivalentes \`a l'existence d'une constante $C\in L$ ind\'ependante de $\mu$ telle que, pour tout $a\in \Q$, tout $j\in \{0,\cdots,k-2\}$ et tout $n\in {\mathbb Z}$ :  
\begin{eqnarray}\label{chaud1} 
\int_{a+p^n\Z}(z-a)^jd\mu(z)&\in & Cp^{n(j-{\rm val}(\alpha))}\O\\  
\label{chaud2}\int_{\Q-(a+p^n\Z)}\left(\frac{\alpha p}{\beta}\right)^{\!\!\!{\rm val}(z-a)}\!\!\!\!\!\!\!(z-a)^{k-2-j}d\mu(z)&\in & Cp^{n({\rm val}(\alpha)-j)}\O
\end{eqnarray}
(si $g=$ {\tiny$\begin{pmatrix}1&\lambda \\0&\mu\end{pmatrix}$}, poser $n=-{\rm val}(\mu)$ et $a=\lambda/\mu$; en fait, on peut prendre $C=1$). De (\ref{chaud1}), on d\'eduit facilement en notant que ${\rm val}(z-a)=n-1$ si $z\in (a+p^{n-1}\Z)-(a+p^n\Z)$ et quitte \`a modifier $C$ :
\begin{equation}\label{chaud3} 
\int_{(a+p^{n-1}\Z)-(a+p^n\Z)}\!\!\!\left(\frac{\alpha p}{\beta}\right)^{{\rm val}(z-a)}(z-a)^{k-2-j}d\mu(z)\in Cp^{n({\rm val}(\alpha)-j)}\O.
\end{equation}
En d\'ecompo\-sant $\Q-(a+p^n\Z)=\Q-(a+p^{n+1}\Z)\setminus (a+p^{n}\Z)-(a+p^{n+1}\Z)$, puis $\Q-(a+p^{n+1}\Z)=\Q-(a+p^{n+2}\Z)\setminus (a+p^{n+1}\Z)-(a+p^{n+2}\Z)$ etc. jusqu'\`a arriver \`a $\Q-(a+p^{n+m}\Z)$ avec $n+m\geq 0$, on d\'eduit de (\ref{chaud3}) que (\ref{chaud2}) pour $j<{\rm val}(\alpha)$ d\'ecoule de (\ref{chaud1}) et de (\ref{chaud2}) pour $n\geq 0$ (utiliser (\ref{chaud3}) pour les morceaux compacts dans la d\'ecomposition et (\ref{chaud2}) avec $n'=n+m\geq 0$ pour le restant). Si $a\ne 0$, en d\'ecompo\-sant $\Q-(a+p^n\Z)=\Q-(a+p^{n-1}\Z)\amalg (a+p^{n-1}\Z)-(a+p^n\Z)$, puis $\Q-(a+p^{n-1}\Z)=\Q-(a+p^{n-2}\Z)\amalg (a+p^{n-2}\Z)-(a+p^{n-1}\Z)$ etc. jusqu'\`a arriver \`a $\Q-(a+p^{n-m}\Z)$ avec $n-m<{\rm val}(a)$ et $n\leq m$, on d\'eduit de (\ref{chaud3}) que (\ref{chaud2}) pour $a\ne 0$ et $j\geq {\rm val}(\alpha)$ d\'ecoule de (\ref{chaud1}) et de (\ref{chaud2}) pour $a=0$ et $n\leq 0$ (utiliser (\ref{chaud3}) pour les morceaux compacts dans la d\'ecomposition puis d\'evelopper $(z-a)^{k-2-j}$ et utiliser (\ref{chaud2}) avec $a=0$ et $n'=n-m\leq 0$ pour le restant). Autrement dit, $(\ref{chaud1})$ et $(\ref{chaud2})$ sont \'equivalents \`a :
\begin{itemize}
\item[(i)] (\ref{chaud1});
\item[(ii)] (\ref{chaud2}) pour $n\geq 0$;
\item[(iii)] (\ref{chaud2}) pour $a=0$ et $n\leq 0$.
\end{itemize}
Par ailleurs, tout $\mu\in \pi(\alpha)^*$ s'\'ecrit $\mu=(\mu_1,\mu_2)$ o\`u $\mu_i\in {\rm Pol}^{k-2}(\Z,L)^*$ (si $f\in \pi(\alpha)$,
$\int_{\Q}f(z)d\mu(z)=\int_{\Z}f_1(z)d\mu_1(z)+\int_{\Z}f_2(z)d\mu_2(z)$). Un calcul facile (laiss\'e au lecteur) montre que $\mu_1$ et $\mu_2$ v\'erifient (\ref{borneamice}) pour $r={\rm val}(\alpha)$ et $d=k-2$ avec $\|\mu_i\|_{{\rm val}(\alpha),k-2}\leq C$ (i.e. $\mu_1, \mu_2\in {\mathcal C}^{{\rm val}(\alpha)}(\Z,L)$ par le th\'eor\`eme \ref{amice2} avec leurs normes born\'ees, i.e. $\mu$ est dans une boule de $B(\alpha)^*\subset \pi(\alpha)^*$) si et seulement si $\mu$ v\'erifie (quitte \`a modifier $C$) :
\begin{equation}\label{ch0} 
\int_{a+p^{n}\Z}(z-a)^jd\mu(z)\in Cp^{n(j-{\rm val}(\alpha))}\O 
\end{equation} 
pour tout $a\in p\Z$, tout $j\in\{0,\cdots,k-2\}$ et tout entier $n\geq 1$, puis :
\begin{equation}\label{ch1} 
\int_{a^{-1}+p^{n-2{\rm val}(a)}\Z}\Big(\frac{p\alpha}{\beta}\Big)^{\!\!{\rm val}
(z)}z^{k-2-j}(1-az)^jd\mu(z)\in Cp^{n(j-{\rm val}(\alpha))}\O 
\end{equation} 
pour tout $a\in \Z-\{0\}$, tout $j\in\{0,\cdots,k-2\}$ et tout entier $n>{\rm val}(a)$, et enfin : 
\begin{equation}\label{ch2} 
\int_{\Q-p^{n}\Z}\Big(\frac{p\alpha}{\beta}\Big)^{\!\!{\rm val}(z)}z^{k-2-j}d\mu(z)\in Cp^{n({\rm val}(\alpha)-j)}\O  
\end{equation} 
pour tout $j\in\{0,\cdots,k-2\}$ et tout entier $n\leq 0$. En d\'eveloppant $z^{k-2-j}=((z-a^{-1})+a^{-1})^{k-2-j}$ dans (\ref{ch1}), un calcul montre que, quitte \`a modifier $C$, (\ref{ch0}), (\ref{ch1}) et (\ref{ch2}) sont \'equivalents \`a :
\begin{itemize}
\item[(iv)] (\ref{chaud1}) pour $a\ne 0$;
\item[(v)] (\ref{chaud1}) pour $a=0$ et $n\geq 0$;
\item[(vi)] (\ref{chaud2}) pour $a=0$ et $n\leq 0$.
\end{itemize}
Si $\mu$ est comme en (\ref{thedual}), i.e. si $\mu$ v\'erifie (i) \`a
(iii), alors {\it a fortiori} $\mu$ v\'erifie (iv) \`a (vi) et donc
$\mu\in B(\alpha)^*\subset \pi(\alpha)^*$. Mais on a plus. En faisant
tendre $n$ vers $-\infty$ dans (\ref{chaud1}) lorsque $a=0$, on voit
que (\ref{chaud1}) pour $j<{\rm val}(\alpha)$ et $a=0$ implique que
$\mu$ annule les fonctions $z^j\in B(\alpha)$. En faisant tendre $n$
vers $+\infty$ dans (\ref{chaud2}), on voit que (\ref{chaud2}) pour
$j<{\rm val}(\alpha)$ implique que $\mu$ annule les fonctions
$\big(\frac{\alpha p}{\beta}\big)^{{\rm val}(z-a)}(z-a)^{k-2-j}\in
B(\alpha)$. Un examen plus approfondi (sans difficult\'e mais que nous
omettons pour ne pas allonger la preuve) montre que les conditions (i)
\`a (iii) pr\'ec\'edentes sont en fait {\it \'equivalentes} aux
conditions (iv) \`a (vi) avec les deux conditions suppl\'ementaires
que $\mu$ annule les fonctions $z^j$ pour $j<{\rm val}(\alpha)$ et les
fonctions $\big(\frac{\alpha p}{\beta}\big)^{{\rm
    val}(z-a)}(z-a)^{k-2-j}$ pour $a\in \Q$ et $j<{\rm val}(\alpha)$,
c'est-\`a-dire les fonctions de $L(\alpha)$. Autrement dit, on obtient
que le Banach dual du compl\'et\'e cherch\'e est isomorphe dans
$\pi(\alpha)^*$ au sous-espace de Banach de $B(\alpha)^*$ form\'e des
$\mu$ qui annulent $L(\alpha)$, c'est-\`a-dire \`a
$(B(\alpha)/L(\alpha))^*$. En particulier, $(B(\alpha)/L(\alpha))^*$
est un $\G$-Banach unitaire. Comme les Banach ne sont pas r\'eflexifs,
nous allons devoir faire un passage par les topologies faibles pour
d\'eduire l'isomorphisme de l'\'enonc\'e. L'injection
$B(\alpha)/L(\alpha)\hookrightarrow ((B(\alpha)/L(\alpha))^*)^*$
\'etant une immersion ferm\'ee $\G$-\'equivariante,
$B(\alpha)/L(\alpha)$ est aussi un $\G$-Banach unitaire. Cela
entra\^\i ne facilement que l'application $\pi(\alpha)\rightarrow
B(\alpha)/L(\alpha)$ induit un morphisme $\G$-\'equivariant continu du
compl\'et\'e unitaire ci-dessus de $\pi(\alpha)$ vers
$B(\alpha)/L(\alpha)$, donc un morphisme continu sur les duaux munis
de leur topologie faible (qui sont des {\og modules compacts \`a
isog\'enie pr\`es \fg} au sens de \cite{ST3}). Mais on vient de voir que ce morphisme sur les duaux \'etait bijectif (et m\^eme un isomorphisme topologique pour les topologies fortes). Par \cite[lemme 4.2.2]{Br2}, on en d\'eduit que c'est aussi un isomorphisme topologique pour les topologies faibles. Par dualit\'e (cf. \cite[th\'eor\`eme 1.2]{ST3}), on obtient l'isomorphisme topologique $\G$-\'equivariant de l'\'enonc\'e. Le cas $\beta$ se traite de m\^eme.
\end{proof}

Rappelons qu'il existe, \`a multiplication par un scalaire non nul pr\`es, un unique morphisme non nul $\G$-\'equivariant :
\begin{equation}\label{entrelisse}
I^{\rm lisse}:{\rm Ind}_{\B}^{\G}{\rm nr}(\beta^{-1})\otimes {\rm nr}(p\alpha^{-1})\rightarrow {\rm Ind}_{\B}^{\G}{\rm nr}(\alpha^{-1})\otimes {\rm nr}(p\beta^{-1})
\end{equation}
qui est un isomorphisme lorsque $\alpha\ne p\beta$ et qui a un noyau et un conoyau de dimension $1$ lorsque $\alpha=p\beta$ (voir \cite[\S4.5]{Bu} par exemple). En termes de fonctions localement constantes sur $\Q$, ce morphisme est donn\'e explicitement par :
\begin{eqnarray}\label{integrale} 
I^{\rm lisse}(h)(z)=\int_{\Q}\!\!\Big(\frac{p\beta}{\alpha}\Big)^{{\rm
 val}(x)}\!h(z+x^{-1})dx &=&\int_{\Q}\!\!\Big(\frac{p\alpha}{\beta}\Big)^{{\rm val}(x)}\!h(z+x)dx\\
\nonumber &=&\int_{\Q}\!\!\Big(\frac{p\alpha}{\beta}\Big)^{{\rm val}(z-x)}\!h(x)dx
\end{eqnarray} 
o\`u $dx$ est la mesure de Haar sur $\Q$ (\`a valeurs dans ${\mathbb Q}\subset L$). Comme la th\'eorie des repr\'esentations lisses est alg\'ebrique, il n'y a pas de probl\`e\-mes de convergence dans les int\'egrales ci-dessus car on peut toujours remplacer les sommes infinies aux voisinages de $0$ ou de $-\infty$ par des expressions alg\'ebriques en $\alpha p\beta^{-1}$ parfaitement d\'efinies. En tensorisant par l'application identit\'e sur ${\rm Sym}^{k-2}L^2$, on en d\'eduit un morphisme non nul $\G$-\'equivariant :
\begin{equation}\label{intertw}
I:\pi(\beta)\longrightarrow \pi(\alpha)
\end{equation}
qui est un isomorphisme lorsque $\alpha\ne p\beta$. 

\begin{coro}\label{entrepadique}
Les repr\'esentations $B(\alpha)/L(\alpha)$ et
$B(\beta)/L(\beta)$ sont des $\G$-Banach unitaires et on a un diagramme commutatif $\G$-\'equiva\-riant :
$$\begin{matrix}B(\beta)/L(\beta)&\buildrel \widehat I\over\longrightarrow &B(\alpha)/L(\alpha)\\
\uparrow &&\uparrow \\ \pi(\beta)&\buildrel I\over\longrightarrow &\pi(\alpha)\end{matrix}$$  
o\`u $I$ est le morphisme $\G$-\'equivariant de (\ref{intertw}). Lorsque $\alpha\ne p\beta$, les fl\`eches $I$ et $\widehat I$ sont des isomorphismes.
\end{coro} 
\begin{proof} 
Cela d\'ecoule du th\'eor\`eme \ref{complete} car l'image par une fl\`eche $\G$-\'equivariante d'un $\O[\G]$-module de type fini est aussi un $\O[\G]$-module de type fini.
\end{proof} 

Lorsque $\alpha=p\beta$ (ce qui implique ${\rm val}(\beta)=(k-2)/2$), il est \'evident que $B(\beta)/L(\beta)$ est non nul puisque $L(\beta)$ est dans ce cas une repr\'esentation de dimension finie (isomorphe via (\ref{algebraic}) \`a la repr\'esentation alg\'ebrique ${\rm Sym}^{k-2}L^2$ \`a une torsion non ramifi\'ee pr\`es). Dans les autres cas, le th\'eor\`eme \ref{complete} ne d\'emontre en rien que les espaces de Banach $B(\alpha)/L(\alpha)$ et $B(\beta)/L(\beta)$ sont non nuls. Mais on a :

\begin{prop}\label{existencereseau}
Si $\alpha\ne p\beta$, le Banach $B(\alpha)/L(\alpha)$ (resp. $B(\beta)/L(\beta)$) est non nul si et seulement si $\pi(\alpha)$ (resp. $\pi(\beta)$) poss\`ede un $\O$-r\'eseau stable par $\G$. Si $\alpha=p\beta$, le Banach $B(\alpha)/L(\alpha)$ est non nul si et seulement si $\pi(\alpha)$ poss\`ede un $\O$-r\'eseau stable par $\G$. 
\end{prop}
\begin{proof}
Rappelons qu'un $\O$-r\'eseau est par d\'efinition un sous-$\O$-module g\'en\'erateur qui ne contient pas de $L$-droite. Supposons d'abord $\alpha\ne p\beta$, de sorte que les repr\'esentations $\pi(\alpha)$ et $\pi(\beta)$ sont (alg\'ebriquement) irr\'eductibles. Si $B(\alpha)/L(\alpha)\ne 0$, l'application canonique $\pi(\alpha)\rightarrow B(\alpha)/L(\alpha)$ est injective car non nulle (car d'image dense) et une boule unit\'e de $B(\alpha)/L(\alpha)$ stable par $\G$ induit un r\'eseau stable par $\G$ sur $\pi(\alpha)$. Inversement, supposons que $\pi(\alpha)$ poss\`ede un $\O$-r\'eseau stable par $\G$, alors pour tout $f$ non nul dans $\pi(\alpha)$, $\O[\G]f\subset \pi(\alpha)$ est un $\O$-r\'eseau de $\pi(\alpha)$ de type fini comme $\O[\G]$-module. Il est g\'en\'erateur car $\pi(\alpha)$ est irr\'eductible et il ne contient pas de $\O$-droite car, \`a multiplication pr\`es par un scalaire, il est contenu dans un $\O$-r\'eseau stable par $\G$ de $\pi(\alpha)$. L'application de $\pi(\alpha)$ dans son compl\'et\'e par rapport \`a $\O[\G]f$, qui est $B(\alpha)/L(\alpha)$ par le th\'eor\`eme \ref{complete}, est  alors injective et en particulier $B(\alpha)/L(\alpha)\ne 0$. Lorsque $\alpha=p\beta$, ce qui suppose $k>2$, $\pi(\alpha)$ n'est plus irr\'eductible et a un quotient isomorphe \`a ${\rm Sym}^{k-2}L^2$. Si $B(\alpha)/L(\alpha)\ne 0$, l'application non nulle $\pi(\alpha)\rightarrow B(\alpha)/L(\alpha)$ reste injective sinon elle induirait une injection non nulle ${\rm Sym}^{k-2}L^2\hookrightarrow B(\alpha)/L(\alpha)$ ce qui est impossible car, pour $k>2$, ${\rm Sym}^{k-2}L^2$ ne poss\`ede pas de $\O$-r\'eseau stable par $\G$.
\end{proof}

Notons que, lorsque $\alpha=p\beta$, $\pi(\beta)$ ne peut poss\'eder de $\O$-r\'eseau stable par $\G$ car sa sous-repr\'esentation irr\'eductible ${\rm Sym}^{k-2}L^2$ n'en poss\`ede pas. Dans ce cas, l'application $\pi(\beta)\rightarrow B(\beta)/L(\beta)$ est non injective (son noyau est pr\'ecis\'ement ${\rm Sym}^{k-2}L^2$). Via la proposition \ref{existencereseau}, on peut montrer pour des petites valeurs de $k$ ou pour les valeurs de $(k,\alpha,\beta)$ provenant des formes modulaires que les Banach $B(\alpha)/L(\alpha)$ et $B(\beta)/L(\beta)$ sont non nuls, voir par exemple \cite{Br1}, \cite[\S1.3]{Br2}, \cite{Br3}, \cite{Em}. On va voir dans la suite que la non nullit\'e pour tout $k$ (et tout $\alpha$, $\beta$) d\'ecoule de la th\'eorie des $(\varphi,\Gamma)$-modules. Une autre approche possible, purement en termes de th\'eorie des repr\'esentations, est pr\'esent\'ee dans \cite[\S2]{Em}.

\section{Repr\'esentations de $\G$ et $(\varphi,\Gamma)$-modules}

Le but de cette partie est de d\'emontrer l'existence d'un isomorphisme topologique canonique entre $(\varprojlim_{\psi} D(V))^{\rm b}$ et le dual $\Pi(V)^*$ (muni de sa topologie faible) et d'en d\'eduire que $\Pi(V)$ est toujours non nul, topologiquement irr\'eductible et admissible. Ces \'enonc\'es \'etaient conjectur\'es (et des cas particuliers d\'emontr\'es) dans \cite{Br1} et \cite{Br2}. Le fait remarquable est que ces \'enonc\'es, enti\`erement du c\^ot\'e ${\rm GL}_2$, se d\'emontrent en passant par le c\^ot\'e {\it galoisien}. On fixe une fois pour toutes une repr\'esentation cristalline irr\'eductible $V$ comme au \S\ref{cristallin} avec $\dcris(V)=D(\alpha,\beta)$.

\subsection{Deux lemmes}

Le but de ce paragraphe est de d\'emontrer deux lemmes techniques mais importants utilis\'es dans les paragraphes suivants. On utilise sans commentaire certaines notations du \S\ref{rcm}.

\begin{lemm}\label{explicit} 
Soit $m\in {\mathbb N}-\{0\}$, $w_{\alpha}$, $w_{\beta}\in \r^+$ et $\mu_{\alpha}$, $\mu_{\beta}$ les distributions localement analytiques sur $\Z$ correspondantes par (\ref{Yvette}). La condition :  
$$\varphi^{-m}(w_{\alpha}\otimes e_{\alpha}+w_{\beta}\otimes
e_{\beta})\in {\rm Fil}^0(L_m[[t]]\otimes_LD(\alpha,\beta))$$  
est \'equivalente aux \'egalit\'es dans $\Qpbar$ :  
$$\alpha^m\int_{\Z}z^j\zeta_{p^m}^zd\mu_{\alpha}(z)
=\beta^m\int_{\Z}z^j\zeta_{p^m}^zd\mu_{\beta}(z)$$  
pour tout $j\in \{0,\cdots,k-2\}$ et toute racine primitive $p^{m}$-i\`eme $\zeta_{p^m}$ de $1$. 
\end{lemm} 
\begin{proof}
On a :
\begin{equation*}
\varphi^{-m}(X)=\zeta_{p^m}{\rm exp}(t/p^m)-1=\zeta_{p^m}({\rm exp}(t/p^m)-1)+\zeta_{p^m}-1
\end{equation*} 
dans $\Q(\pmb{\mu}_{p^m})[[t]]$ (voir \S\ref{rcm}). En posant $w_{\alpha}=\sum_{i=0}^{+\infty}\alpha_iX^i$ et
$w_{\beta}=\sum_{i=0}^{+\infty}\beta_iX^i$ ($\alpha_i$, $\beta_i\in L$), la condition sur ${\rm Fil}^0$ est \'equivalente \`a :  
\begin{multline} \label{filo}
\alpha^m\sum_{i=0}^{+\infty}\alpha_i\varphi^{-m}(X)^i 
e_{\alpha}+\beta^m\sum_{i=0}^{+\infty}\beta_i\varphi^{-m}(X)^ie_{\beta}\in\\
L_m[[t]](e_{\alpha}+e_{\beta})\oplus   
({\rm exp}(t/p^m)-1)^{k-1}(L_m[[t]]e_{\alpha}) 
\end{multline} 
en notant que ${\rm exp}(t/p^m)-1$ engendre ${\rm gr}^1(\Q[[t]])=\Q\overline
t$. On peut supposer $L$ aussi grand que l'on veut en (\ref{filo}), et en particulier contenant $\Q(\pmb{\mu}_{p^m})$. En utilisant :
$$L_m=\Q(\pmb{\mu}_{p^m})\otimes_{\Q}L=\prod_{\Q(\pmb{\mu}_{p^m})\hookrightarrow L} L$$ 
et en d\'eveloppant $\varphi^{-m}(X)^i=(\zeta_{p^m}({\rm exp}(t/p^m)-1)+\zeta_{p^m}-1)^i$, 
un calcul facile montre que la condition (\ref{filo}) est \'equivalente aux \'egalit\'es dans $L$ :
\begin{equation}\label{combin} 
\alpha^m\sum_{i=j}^{+\infty}\alpha_i\binom{i}{j}(\zeta_{p^m}-1)^{i-j}
=\beta^m\sum_{i=j}^{+\infty}\beta_i\binom{i}{j}(\zeta_{p^m}-1)^{i-j} 
\end{equation} 
pour tout $j\in \{0,\cdots,k-2\}$ et {\it toute} racine primitive $p^m$-i\`eme $\zeta_{p^m}$ de $1$. Noter que les s\'eries en (\ref{combin}) convergent bien car ${\rm val}(\zeta_{p^m}-1)>0$. En se souvenant que $\alpha_i=\int_{\Z}\binom{z}{i}d\mu_{\alpha}(z)$ et en utilisant le d\'eveloppement de Mahler (\ref{mahler}) :  
$$\binom{z}{j}\zeta_{p^m}^{z-j}=\sum_{i=j}^{+\infty}\binom{i}{j}(\zeta_{p^m}-1)^{i-j}\binom{z}{i}$$  
on obtient : 
\begin{eqnarray*} 
\sum_{i=j}^{+\infty}\alpha_i\binom{i}{j}(\zeta_{p^m}-1)^{i-j}&=&\int_{\Z}
\Big(\sum_{i=j}^{+\infty}\binom{z}{i}\binom{i}{j}(\zeta_{p^m}-1)^{i-j}\Big)d\mu_{\alpha}(z)\\  
&=&\zeta_{p^m}^{-j}\int_{\Z}\binom{z}{j}\zeta_{p^m}^{z}d\mu_{\alpha}(z) 
\end{eqnarray*} 
(la s\'erie $\sum_{i=j}^{n}\binom{i}{j}(\zeta_{p^m}-1)^{i-j}\binom{z}{i}$ convergeant vers $\sum_{i=j}^{+\infty}\binom{i}{j}(\zeta_{p^m}-1)^{i-j}\binom{z}{i}$ dans ${\rm An}(\Z,L)$ (cf. \cite[\S2.1.2]{Co2}), on peut inverser $\int$ et $\sum$). On a la m\^eme \'egalit\'e avec $\beta_i$ et $\mu_{\beta}$. Avec (\ref{combin}), on d\'eduit le r\'esultat. 
\end{proof} 

Via l'identification $\Q/\Z={\mathbb Z}[1/p]/{\mathbb Z}$, on peut d\'efinir le nombre complexe alg\'ebrique $e^{2i\pi z}$ pour tout $z\in \Q$ (par exemple, $e^{2i\pi z}=1$ si $z\in \Z$). En fixant des plongements $\overline{\mathbb Q}\hookrightarrow {\mathbb C}$ et $\overline{\mathbb Q}\hookrightarrow \Qpbar$, on peut voir $e^{2i\pi z}$ comme un \'el\'ement de $\Qpbar$. On obtient ainsi un caract\`ere additif localement constant $\Q\rightarrow \Qpbar^{\times}$, $z\mapsto e^{2i\pi z}$ trivial sur $\Z$. Ce que l'on fera dans la suite ne d\'ependra pas du choix de ce caract\`ere, i.e. du choix des plongements.

Notons ${\rm Pol}^d(\Q,L)$ le $L$-espace vectoriel des fonctions localement polynomiales \`a support compact $f:\Q\rightarrow L$ de degr\'e (local) au plus $d$. Si $\mu$ est une forme lin\'eaire sur ${\rm Pol}^d(\Q,L)$, $U$ un ouvert compact de $\Q$ et $f:\Q\rightarrow L$ une fonction localement polynomiale de degr\'e au plus $d$ \`a support quelconque, on note comme d'habitude $\int_Uf(z)d\mu(z)\=\mu({\bf 1}_U(z)f(z))$.

\begin{lemm}\label{crucial} 
Soit $\mu_{\alpha}$ et $\mu_{\beta}$ deux formes lin\'eaires sur ${\rm Pol}^{k-2}(\Q,L)$. Les
\'enonc\'es suivants sont \'equivalents :
\begin{itemize}
\item[(i)] Pour tout $j\in \{0,\cdots,k-2\}$, tout $y\in \Q^{\times}$ et tout $N>{\rm val}(y)$, on a : 
\begin{equation}\label{rajout} 
\int_{p^{-N}\Z}z^je^{2i\pi zy}d\mu_{\beta}(z)
=\Big(\frac{\beta}{\alpha}\Big)^{{\rm val}(y)}\int_{p^{-N}\Z}z^je^{2i\pi zy}d\mu_{\alpha}(z). 
\end{equation} 
\item[(ii)] Pour tout $f\in \pi(\alpha)$ \`a support compact (comme fonction sur $\Q$ via (\ref{commefonction})) tel que $I(f)\in \pi(\beta)$ est aussi \`a support compact (comme fonction sur $\Q$ via (\ref{commefonction})), on a :  
\begin{equation}\label{debutinter} 
\int_{\Q}f(z)d\mu_{\beta}(z)=\frac{1-\frac{\alpha}{\beta}}{1-\frac{\beta}{p\alpha}}\int_{\Q}I(f)(z)d\mu_{\alpha}(z).  
\end{equation}
\end{itemize}
\end{lemm} 
\begin{proof}
Comme $\alpha\ne \beta$ et $\beta\ne p\alpha$, la constante dans (\ref{debutinter}) est bien d\'efinie et non nulle. Soit $h:\Q\rightarrow L$ une fonction localement constante \`a support compact et $\widehat h$ la transform\'ee de Fourier (usuelle) de $h$. Rappelons que $\widehat h$ est aussi une fonction localement constante sur $\Q$ \`a support compact telle que $\widehat h(x)=\int_{\Q}h(z)e^{-2i\pi zx}dz$ et $h(z)=\int_{\Q}\widehat h(x)e^{2i\pi zx}dx$ o\`u $dx$, $dz$ d\'esignent la mesure de Haar sur $\Q$. Pour $|z|\gg 0$, on a $I^{\rm lisse}(h)(z)=(\alpha p\beta^{-1})^{{\rm val}(z)}\int_{\Q}h(x)dx=(\alpha p\beta^{-1})^{{\rm val}(z)}\widehat h(0)$ et on voit que $I^{\rm lisse}(h)$ est \`a support compact dans $\Q$ si et seulement si $\widehat h(0)=0$. Supposons donc $\widehat h(0)=0$ et soit $N\in {\mathbb N}$ tel que $h$ et $I^{\rm lisse}(h)$ ont leur support dans $p^{-N}\Z$ et tel que $\widehat h|_{p^N\Z}=0$. Pour $j\in \{0,\cdots,k-2\}$ et $z\in p^{-N}\Z$, on a :  
\begin{eqnarray*} 
I(z^jh)(z)&=&z^jI^{\rm lisse}(h)(z)\ =\ z^j\int_{p^{-N}\Z}(\alpha p\beta^{-1})^{{\rm val}(x)}h(z+x)dx\\ 
&=&z^j\int_{p^{-N}\Z}(\alpha p\beta^{-1})^{{\rm val}(x)}
\bigg(\int_{\Q-p^N\Z}\widehat h(y)e^{2i\pi y(z+x)}dy\bigg)dx\\ 
&=&z^j\int_{\Q-p^N\Z}\widehat h(y)e^{2i\pi zy}\bigg(\int_{p^{-N}\Z}
(\alpha p\beta^{-1})^{{\rm val}(x)}e^{2i\pi xy}dx\bigg)dy. 
\end{eqnarray*} 
En d\'ecomposant $p^{-N}\Z=p^{-N}\Z^{\times}\amalg p^{-N+1}
\Z^{\times}\amalg\cdots$, un calcul facile (laiss\'e au lecteur) donne pour $N>{\rm val}(y)$ : 
$$\int_{p^{-N}\Z}\!\!\!\!\!(\alpha p\beta^{-1})^{{\rm val}(x)}e^{2i\pi xy}dx=\!\!\sum_{i=-N}^{+\infty}(\alpha
p\beta^{\!-1})^i\!\!\!\int_{p^i\Z^{\times}}\!\!\!e^{2i\pi xy}dx=\frac{1-\frac{\beta}{p\alpha}}{1-\frac{\alpha}
{\beta}}\Big(\frac{\beta}{\alpha}\Big)^{{\rm val}(y)}$$  
d'o\`u : 
\begin{equation}\label{fourier1} 
I(z^jh)(z)=\frac{1-\frac{\beta}{p\alpha}}{1-\frac{\alpha}
{\beta}}\int_{\Q-p^N\Z}\widehat h(y)z^je^{2i\pi zy}\Big(\frac{\beta}
{\alpha}\Big)^{{\rm val}(y)}dy 
\end{equation} 
pour $z\in p^{-N}\Z$ et $I(z^jh)(z)=0$ sinon. De m\^eme, on a :  
\begin{equation}\label{fourier2} 
z^jh(z)=\int_{\Q-p^N\Z}\widehat h(y)z^je^{2i\pi zy}dy 
\end{equation} 
pour $z\in p^{-N}\Z$ et $z^jh(z)=0$ sinon. Notons que (\ref{fourier1}) and (\ref{fourier2}) sont en fait des sommes finies sur le m\^eme ensemble (fini) de valeurs de $y$. En rempla\c cant $I(z^jh)(z)$ et $z^jh(z)$ dans (\ref{variant}) ci-dessous par les sommes finies (\ref{fourier1}) et (\ref{fourier2}), on voit que (i) entra\^\i ne (ii). R\'eciproquement, supposons que pour tout $j\in \{0,\cdots,k-2\}$, tout $N\in \N$ et toute fonction $h$ comme ci-dessus localement constante \`a support dans $p^{-N}\Z$, on a :
 \begin{equation}\label{variant} 
\int_{p^{-N}\Z}I(z^jh)(z)d\mu_{\alpha}(z)=\frac{1-\frac{\beta}
{p\alpha}}{1-\frac{\alpha}{\beta}}\int_{p^{-N}\Z}z^jh(z)d\mu_{\beta}(z).
\end{equation} 
Soit $y\in \Q^{\times}$ tel que $N>{\rm val}(y)$, $\widehat h(z)\= {\bf 1}_{y+p^N\Z}(z)$ et :
$$h(z)\=\int_{\Q}\widehat h(x)e^{2i\pi zx}dx=\int_{p^{N}\Z}e^{2i\pi z(y+x)}dx=\frac{1}{p^N}{\bf 1}_{p^{-N}\Z}(z)e^{2i\pi zy}.$$
Alors $I(z^jh)$ est aussi \`a support compact car $\widehat h|_{p^N\Z}=0$ et un calcul via (\ref{fourier1}) montre que :
$$I(z^jh)(z)=\frac{1}{p^N}\frac{1-\frac{\beta}{p\alpha}}{1-\frac{\alpha}
{\beta}}\left(\frac{\beta}{\alpha}\right)^{{\rm val}(y)}\!\!\!{\bf 1}_{p^{-N}\Z}(z)z^je^{2i\pi zy}.$$
On peut donc appliquer l'\'egalit\'e (\ref{variant}) \`a $h$ qui est alors exactement l'\'egalit\'e (\ref{rajout}) multipli\'ee par $p^{-N}$. Cela montre que (ii) entra\^\i ne (i) et ach\`eve la preuve.
\end{proof} 

\subsection{De $\g$ vers $\G$}

Le but de ce paragraphe est de construire une application continue $(\varprojlim_{\psi}D(V))^{\rm b}\rightarrow \Pi(V)^*$.

Soit $T\subset V$ un $\O$-r\'eseau stable par $\g$. On reprend les notations du \S\ref{phigamma}, en particulier on dispose du $\O[[X]]$-module de type fini $D^{\sharp}(T)$ muni de la surjection $\psi:D^{\sharp}(T)\twoheadrightarrow D^{\sharp}(T)$ et de l'action semi-lin\'eaire de $\Gamma$ qui commute \`a $\psi$. On dispose aussi de l'isomorphisme topologique de la proposition \ref{ddiese} qui permet de remplacer $(\varprojlim_{\psi}D(V))^{\rm b}$ par $(\varprojlim_{\psi}D^{\sharp}(T))\otimes_{\O}L$ et on sait par le corollaire \ref{ddiesedcrisproj} que $(\varprojlim_{\psi}D^{\sharp}(T))\otimes_{\O}L$ co\"\i ncide avec les suites d'\'el\'ements $w_{\alpha,n}\otimes e_{\alpha} + w_{\beta,n}\otimes e_{\beta}$ de $\r^+\otimes_L\dcris(V)$ telles que :
\begin{itemize}
\item[(i)] $\forall\ n\geq 0$, $w_{\alpha,n}$ (resp. $w_{\beta,n}$) est d'ordre ${\rm val}(\alpha)$ (resp. ${\rm val}(\beta)$) dans $\r^+$ et $\|w_{\alpha,n}\|_{{\rm val}(\alpha)}$ (resp. $\|w_{\beta,n}\|_{{\rm val}(\beta)}$) est born\'e ind\'ependamment de $n$;
\item[(ii)] $\forall\ n\geq 0$ et $\forall\ m\geq 1$, on a :
$$\varphi^{-m}(w_{\alpha,n}\otimes e_{\alpha}+w_{\beta,n}\otimes 
e_{\beta})\in {\rm Fil}^0(L_m)[[t]]\otimes_L\dcris(V));$$ 
\item[(iii)] $\forall\ n\geq 1$, $\psi(w_{\alpha,n})=\alpha^{-1}w_{\alpha,n-1}$ et  
$\psi(w_{\beta,n})=\beta^{-1}w_{\beta,n-1}$.
\end{itemize}

Nous allons d'abord d\'efinir une application $L$-lin\'eaire $(\varprojlim_{\psi}D^{\sharp}(T))\otimes_{\O}L\rightarrow \pi(\alpha)^*$. Soit $\mu_{\alpha,n}$ et $\mu_{\beta,n}$ les distributions sur $\Z$ correspondant \`a $\alpha^nw_{\alpha,n}$ et $\beta^nw_{\beta,n}$ par la transform\'ee d'Amice-Mahler (\ref{Yvette}). On associe \`a $(\mu_{\alpha,n})_n$ et $(\mu_{\beta,n})_n$ deux distributions localement analytiques $\mu_{\alpha}$ et $\mu_{\beta}$ sur $\Q$ \`a support compact (i.e. deux formes lin\'eaires continues sur l'espace vectoriel des fonctions localement analytiques sur $\Q$ \`a support compact) en posant :
\begin{equation}\label{prolong} 
\int_{U}f(z)d\mu_{\alpha}(z)\=\int_{\Z}{\bf 1}_U(z/p^N)f(z/p^N)d\mu_{\alpha,N}(z)
\end{equation} 
(resp. avec $\beta$ au lieu de $\alpha$) o\`u $f:\Q\rightarrow L$ est localement analytique (\`a support quelconque) et $U$ est un ouvert compact de $\Q$ contenu dans $p^{-N}\Z$.

\begin{lemm}
La valeur $\int_{\Z}{\bf 1}_U(z/p^N)f(z/p^N)d\mu_{\alpha,N}(z)$ ne d\'epend pas du choix de $N$ tel que $U$ est contenu dans $p^{-N}\Z$.
\end{lemm}
\begin{proof}
Si $\mu$ est une distribution localement analytique sur $\Z$ correspondant \`a $w\in \r^+$ par (\ref{Yvette}), il est facile de voir que la distribution localement analytique $\psi(\mu)$ correspondant \`a $\psi(w)$ v\'erifie :
\begin{equation}\label{psidis}
\int_{\Z}f(z)d\psi(\mu)(z)=\int_{p\Z}f(z/p)d\mu(z).
\end{equation}
On a donc :
\begin{eqnarray*}
\int_{\Z}{\bf 1}_U(z/p^N)f(z/p^{N})d\mu_{\alpha,N}(z)&=&\int_{\Z}{\bf 1}_U(z/p^N)f(z/p^{N})d\psi(\mu_{\alpha,N+1})(z)\\
&\buildrel (\ref{psidis})\over =&\int_{p\Z}{\bf 1}_U(z/p^{N+1})f(z/p^{N+1})d\mu_{\alpha,N+1}(z)\\
&=&\int_{\Z}{\bf 1}_U(z/p^{N+1})f(z/p^{N+1})d\mu_{\alpha,N+1}(z),
\end{eqnarray*}
en remarquant que ${\bf 1}_U(z/p^{N+1})f(z/p^{N+1})$ est \`a support dans $p\Z$.
\end{proof}

Par le lemme \ref{explicit}, la condition (ii) pr\'ec\'edente sur $(w_{\alpha,n},w_{\beta,n})_n$ est
\'equivalente aux \'egalit\'es :  
\begin{equation}\label{filtr} 
\alpha^{m-n}\int_{\Z}z^j\zeta_{p^m}^zd\mu_{\alpha,n}(z)=\beta^{m-n}\int_{\Z}z^j\zeta_{p^m}^zd\mu_{\beta,n}(z)  
\end{equation} 
pour tout $j\in \{0,\cdots,k-2\}$, tout $n\geq 0$ et tout $m\geq 1$. 
  
\begin{coro}\label{equival} 
Avec les notations pr\'ec\'edentes, la condition (ii) ci-dessus sur $(w_{\alpha,n},w_{\beta,n})_n$ est \'equivalente aux \'egalit\'es dans $\Qpbar$ : 
$$\int_{p^{-N}\Z}z^je^{2i\pi zy}d\mu_{\alpha}(z)=\Big(\frac{\alpha}{\beta}\Big)^{{\rm val}(y)}\int_{p^{-N}\Z}z^je^{2i\pi zy}d\mu_{\beta}(z)$$  
pour tout $j\in \{0,\cdots,k-2\}$, tout $y\in \Q^{\times}$ et tout $N>{\rm val}(y)$. 
\end{coro} 
\begin{proof} 
Cela r\'esulte de (\ref{prolong}) et de (\ref{filtr}) en remarquant que $e^{2i\pi y/p^N}$ est une racine primitive d'ordre $p^{N-{\rm val}(y)}$ de $1$.  
\end{proof} 
 
Par le lemme \ref{crucial}, on a donc : 
$$\int_{\Q}f(z)d\mu_{\beta}(z)=\frac{1-\frac{\alpha}{\beta}}{1-\frac{\beta}{p\alpha}}\int_{\Q}I(f)(z)d\mu_{\alpha}(z)$$ 
pour $f\in \pi(\alpha)$ \`a support compact tel que $I(f)\in \pi(\beta)$ est aussi \`a support compact. 

\begin{lemm}\label{versfin} 
Il y a une mani\`ere unique de prolonger $\mu_{\beta}$ et $\mu_{\alpha}$ comme
\'el\'ements respectivement de $\pi(\beta)^*$ et $\pi(\alpha)^*$ telle que, pour tout $f\in \pi(\beta)$ (vue comme fonction sur $\Q$ par (\ref{commefonction})) :  
\begin{equation}\label{extension} 
\int_{\Q}f(z)d\mu_{\beta}(z)=\frac{1-\frac{\alpha}{\beta}}{1-\frac{\beta}{p\alpha}}\int_{\Q}I(f)(z)d\mu_{\alpha}(z). 
\end{equation} 
\end{lemm} 
\begin{proof} 
Un calcul facile \`a partir de (\ref{integrale}) donne pour $j\in \{0,\cdots,k-2\}$ : 
\begin{eqnarray}\label{calculc} 
I(z^j{\mathbf 1}_{\Z})&\!\!\!=\!\!\!&\frac{1-\frac{1}{p}}{1-\frac{\alpha}
{\beta}}z^j{\mathbf 1}_{\Z}(z)+z^j\Big(\frac{p\alpha}
{\beta}\Big)^{\!\!{\rm val}(z)}\!\!\!\!(z-{\mathbf 1}_{\Z}(z))\\
\label{calculnc}
I\Big(z^j\Big(\frac{p\beta}{\alpha}\Big)^{\!\!{\rm val}
(z)}\!\!\!\!(z-{\mathbf 1}_{p\Z}(z))\Big)&\!\!\!=\!\!\!&z^j{\mathbf
  1}_{p\Z}(z)+\frac{1-\frac{1}{p}}{1-\frac{\alpha}{\beta}}z^j
\Big(\frac{p\alpha}{\beta}\Big)^{\!\!{\rm val}(z)}\!\!\!\!(z-{\mathbf 1}_{p\Z}(z)). 
\end{eqnarray} 
Comme $\int_{\Q}z^j{\mathbf 1}_{\Z}d\mu_{\beta}(z)$ et
$\int_{\Q}z^j{\mathbf 1}_{\Z}d\mu_{\alpha}(z)$ sont bien d\'efinis, on voit
en utilisant (\ref{extension}) et (\ref{calculc}) que :
$$\int_{\Q}z^j\Big(\!\frac{p\alpha}{\beta}\!\Big)^{\!\!{\rm val}(z)}\!(z-{\mathbf 1}_{\Z}(z))d\mu_{\alpha}(z)$$
est uniquement d\'etermin\'e, puis en utilisant (\ref{calculnc}) (et (\ref{extension})) que :
$$\int_{\Q}z^j\Big(\!\frac{p\beta}{\alpha}\!\Big)^{\!\!{\rm val}(z)}\!(z-{\mathbf 1}_{p\Z}(z))d\mu_{\beta}(z)$$
est aussi uniquement d\'etermin\'e. Il est facile de v\'erifier que cela d\'efinit bien un prolongement unique de $\mu_{\beta}$ et $\mu_{\alpha}$ comme formes lin\'eaires sur respectivement $\pi(\beta)$ et $\pi(\alpha)$. \end{proof} 

Lorsque $\alpha=p\beta$, on peut voir que la distribution $\mu_{\beta}\in \pi(\beta)^*$ du lemme \ref{versfin} est nulle contre ${\rm Sym}^{k-2}L^2\subset \pi(\beta)$. 

Avec les notations pr\'ec\'edentes, on d\'eduit du lemme \ref{versfin} une application $L$-lin\'eaire :
\begin{eqnarray}\label{versbanach}
(\varprojlim_{\psi}D^{\sharp}(T))\otimes_{\O}L&\longrightarrow &\pi(\alpha)^*\\
\nonumber (w_{\alpha,n}\otimes e_{\alpha}+w_{\beta,n}\otimes
e_{\beta})_n&\longmapsto & \mu_{\alpha}\ {\rm prolong\acute e}
\end{eqnarray}

\begin{lemm}\label{actiongl2}
Soit $\gamma=[a]\in \Gamma\simeq \Z^{\times}$ (cf. \S\ref{borlim}), $z\in \Z$, $(v_n)_n\in \varprojlim_{\psi}D^{\sharp}(T)$ et $\mu_{\alpha}\in \pi(\alpha)^*$ l'image de $(v_n)_n$ par (\ref{versbanach}). Alors :
\begin{itemize}
\item[(i)] $(\psi(v_n))_n$ s'envoie sur $\tiny{\begin{pmatrix}1&0\\0&p\end{pmatrix}}\cdot \mu_{\alpha}$;
\item[(ii)] $(\gamma (v_n))_n$ s'envoie sur $\tiny{\begin{pmatrix}1&0\\0&a\end{pmatrix}}\cdot \mu_{\alpha}$;
\item[(iii)] $(\varphi^n((1+X)^z)v_n)_n$ s'envoie sur $\tiny{\begin{pmatrix}1&z\\0&1\end{pmatrix}}\cdot \mu_{\alpha}$.
\end{itemize}
\end{lemm}
\begin{proof}
Cela d\'ecoule de (\ref{prolong}) et de propri\'et\'es simples de la transform\'ee d'Ami\-ce-Mahler (voir par exemple \cite[\S2.2.2]{Co2}). Nous laissons les d\'etails en exercice au lecteur.
\end{proof}

En particulier, le lemme \ref{actiongl2} induit une action du
groupe $\B$ sur $\varprojlim_{\psi}D^{\sharp}(T)$, qui co\"{\i}ncide
bien s\^ur avec celle de la d\'efinition \ref{actg} (en faisant agir les scalaires par multiplication par le caract\`ere central de $\pi(\alpha)^*$).

\begin{lemm}\label{unsens}
L'application (\ref{versbanach}) se factorise par une injection continue $\B$-\'equivariante :
$$(\varprojlim_{\psi}D^{\sharp}(T))\otimes_{\O}L\hookrightarrow (B(\alpha)/L(\alpha))^*$$
(pour la topologie faible sur $(B(\alpha)/L(\alpha))^*$).
\end{lemm}
\begin{proof}
L'injectivit\'e d\'ecoule via (\ref{prolong}) de l'injectivit\'e dans l'isomorphisme $\r^+\buildrel\sim\over\rightarrow {\rm An}(\Z,L)^*$ (cf. (\ref{Yvette})) et la $\B$-\'equivariance du lemme \ref{actiongl2}. Montrons que l'application $\varprojlim_{\psi}D^{\sharp}(T)\rightarrow \pi(\alpha)^*$ est continue. Notons $\pi(\alpha)_{\rm c}\subset \pi(\alpha)$ (resp. $\pi(\beta)_{\rm c}\subset \pi(\beta)$) le sous-$L$-espace vectoriel des fonctions $f\in \pi(\alpha)$ (resp. $f\in \pi(\beta)$) \`a support compact dans $\Q$. Par (\ref{calculc}), on voit que la fl\`eche $\pi(\alpha)_{\rm c}\oplus \pi(\beta)_{\rm c}\buildrel {\rm incl} \oplus I\over \longrightarrow \pi(\alpha)$ est surjective et induit une immersion ferm\'ee entre espaces de Fr\'echet :
$$\pi(\alpha)^*\hookrightarrow \pi(\alpha)_{\rm c}^*\oplus \pi(\beta)_{\rm c}^*.$$
Il suffit donc de montrer la continuit\'e des deux applications $\varprojlim_{\psi}D^{\sharp}(T)\rightarrow \pi(\alpha)_{\rm c}^*$ et $\varprojlim_{\psi}D^{\sharp}(T)\rightarrow \pi(\beta)_{\rm c}^*$, ce qui d\'ecoule apr\`es passage \`a la limite projective via (\ref{prolong}) de la continuit\'e de $\r^+\buildrel\sim\over\rightarrow {\rm An}(\Z,L)^*$ et de celle de l'injection $D^{\sharp}(T)\hookrightarrow \r^+\otimes_L\dcris(V)$ (cf. remarques \ref{continj} et \ref{weakispt}). Notons $B(V)$ l'espace de Banach dual du module compact $\varprojlim_{\psi}D^{\sharp}(T)$ par l'anti-\'equivalence de cat\'egorie de \cite[\S1]{ST3}. Il est muni d'une action continue unitaire de $\B$ par la proposition \ref{actgcont} (on peut utiliser les arguments de dualit\'e de la preuve de \cite[proposition 1.6]{ST3} pour la continuit\'e de l'action) et on a par ce qui pr\'ec\`ede un morphisme $\B$-\'equivariant (continu) $\pi(\alpha)\rightarrow B(V)$. Par la propri\'et\'e universelle du compl\'et\'e de $\pi(\alpha)$ par rapport \`a un sous-$\O[\B]$-module g\'en\'erateur de type fini et par le th\'eor\`eme \ref{complete}, ce morphisme s'\'etend par continuit\'e en un morphisme $\B$-\'equivariant continu $B(\alpha)/L(\alpha)\rightarrow B(V)$. En redualisant, ce dernier induit un morphisme continu $(\varprojlim_{\psi}D^{\sharp}(T))\otimes_{\O}L\rightarrow (B(\alpha)/L(\alpha))^*$ qui est le morphisme de l'\'enonc\'e.
\end{proof}

\subsection{De $\G$ vers $\g$}\label{laouiw}

Le but de ce paragraphe est de construire une application continue $\Pi(V)^*\rightarrow (\varprojlim_{\psi} D(V))^{\rm b}$ inverse de la pr\'ec\'edente. 

Soit $\mu_{\alpha}\in (B(\alpha)/L(\alpha))^*$ et $\mu_{\beta}\=\frac{1-\frac{\alpha}{\beta}}{1-\frac{\beta}{p\alpha}}\widehat I\circ \mu_{\alpha}\in (B(\beta)/L(\beta))^*$ o\`u $\widehat I$ est le morphisme $\G$-\'equivariant du corollaire \ref{entrepadique}. On d\'efinit une suite $(\mu_{\alpha,n})_n$ de distributions localement analytiques sur $\Z$ en posant :
\begin{equation}\label{prolong2} 
\int_{\Z}f(z)d\mu_{\alpha,n}(z)\=\int_{p^{-n}\Z}f(p^nz)d\mu_{\alpha}(z)
\end{equation} 
et on d\'efinit de m\^eme $(\mu_{\beta,n})_n$. Soit $w_{\alpha,n},\ w_{\beta,n}\in \r^+$ les \'el\'ements correspondant \`a $\alpha^{-n}\mu_{\alpha,n},\ \beta^{-n}\mu_{\beta,n}$ par (\ref{Yvette}). 

\begin{lemm}\label{routine}
La suite d'\'el\'ements $w_{\alpha,n}\otimes e_{\alpha} + w_{\beta,n}\otimes e_{\beta}$ de $\r^+\otimes_L\dcris(V)$ satisfait les conditions (i), (ii) et (iii) du corollaire \ref{ddiesedcrisproj}.
\end{lemm}
\begin{proof}
La condition (iii) est \'evidente \`a partir de (\ref{psidis}) et la condition (ii) d\'ecoule des d\'efinitions, du lemme \ref{crucial} et du corollaire \ref{equival}. V\'erifions la condition (i). Revenant \`a la preuve du th\'eor\`eme \ref{complete}, on a en particulier que $\mu_{\alpha}$ satisfait (\ref{chaud1}) ce qui entra\^\i ne :
\begin{eqnarray*} 
\alpha^{-N}\int_{a+p^n\Z}(z-a)^jd\mu_{\alpha,N}(z)&=&\alpha^{-N}p^{Nj}
\int_{p^{-N}a+p^{n-N}\Z}(z-p^{-N}a)^jd\mu_{\alpha}(z)\\ 
&\in & C_{\mu_{\alpha}}p^{-N{\rm val}(\alpha)}p^{Nj}p^{(n-N)(j-{\rm val}(\alpha))}\O\\ 
&\in &C_{\mu_{\alpha}}p^{n(j-{\rm val}(\alpha))}\O. 
\end{eqnarray*}
pour tout $a\in \Z$, tout $j\in \{0,\cdots, k-2\}$ et tout $n\in \N$. Avec les notations du \S\ref{analyse}, cela entra\^\i ne pour tout $N\in \N$ :  
$$\|\alpha^{-N}\mu_{\alpha,N}\|_{{\rm val}(\alpha),k-2}\leq
c|C_{\mu_{\alpha}}|$$  
pour une constante $c\in \R_{\geq 0}$. On a une borne analogue pour les $\mu_{\beta,N}$. On en d\'eduit (i). 
\end{proof}

Par le lemme \ref{routine} et le corollaire \ref{ddiesedcrisproj}, on a une application $L$-lin\'eaire :
$$(B(\alpha)/L(\alpha))^*\longrightarrow (\varprojlim_{\psi}D^{\sharp}(T))\otimes_{\O}L$$
et il est imm\'ediat \`a partir des d\'efinitions et du lemme \ref{versfin} de v\'erifier qu'elle est inverse de celle du lemme \ref{unsens}.

\begin{theo}\label{phigammabanach}
Il y a un unique isomorphisme topologique (\`a multiplication pr\`es par un scalaire non nul) entre les $L$-espaces vectoriels localement convexes (pour la topologie faible des deux c\^ot\'es): 
$$\big(\varprojlim_{\psi}D(V)\big)^{\rm b}\buildrel\sim \over\longrightarrow\Pi(V)^*$$ 
tel que l'action de {\tiny$\begin{pmatrix}1&0\\0&p^{\mathbb Z}\end{pmatrix}$} sur $\Pi(V)^*$ correspond \`a $(v_n)_n\mapsto (\psi^{\mathbb Z}(v_n))_n$, l'action de {\tiny$\begin{pmatrix}1&0\\0&\Z^{\times}\end{pmatrix}$} \`a celle de $\Gamma$ et l'action de {\tiny$\begin{pmatrix}1&\Z\\0&1\end{pmatrix}$} \`a $(v_n)_n\mapsto ((1+X)^{p^n\Z}v_n))_n$. 
\end{theo}
\begin{proof}
L'existence d'un tel isomorphisme d\'ecoule des r\'esultats
pr\'ec\'e\-dents, sachant qu'une application bijective continue entre
deux {\og modules compacts \`a isog\'enie pr\`es \fg} est un isomorphisme
topologique (c'est la version duale par \cite{ST3} du th\'eor\`eme de
l'image ouverte entre espaces de Banach). Il reste \`a d\'emontrer
l'unicit\'e (\`a scalaire pr\`es) mais cela r\'esulte de la
proposition \ref{schur}.
\end{proof}

Rappelons que $H^i_{\rm Iw}(\Q,V)\= L \otimes_{\O} 
\varprojlim_n H^i({\rm Gal}(\Qpbar/F_n),T)$ 
o\`u $T$ est un $\O$-r\'eseau quelconque
de $V$ stable par $\g$ (voir la d\'efinition \ref{defiiw}).

\begin{coro}\label{iwasawa}
On a un isomorphisme de $\O[[\Z^{\times}]]$-modules :
$$H^1_{\rm Iw}(\Q,V)\simeq {\Pi(V)^*}^{\tiny\begin{pmatrix}
1&0\\0&p^{\mathbb Z}\end{pmatrix}}$$
o\`u $\Z^{\times}$ agit via l'action de $\Gamma$ \`a gauche et via l'action de {\tiny$\begin{pmatrix}1&0\\0&\Z^{\times}\end{pmatrix}$} \`a droite.
\end{coro}
\begin{proof}
Cela d\'ecoule du th\'eor\`eme \ref{phigammabanach} et de la proposition \ref{cciw}.
\end{proof}

\begin{rema}\label{iwasawa2}
Le $\O[[\Z^{\times}]]$-module $H^2_{\rm Iw}(\Q,V)$ s'identifie aussi
aux coinvariants de $\Pi(V)^*$ sous l'action de
{\tiny$\begin{pmatrix}1&0\\0&p^{\mathbb Z}\end{pmatrix}$} et en fait, ces
deux espaces sont nuls par le (ii) de la proposition
\ref{bpriw} parce que $V$ est irr\'eductible. 
En effet, le th\'eor\`eme \ref{phigammabanach} et le 
corollaire \ref{h2ddiese} nous disent que les
coinvariants d'un r\'eseau de $\Pi(V)^*$ sous l'action de
{\tiny$\begin{pmatrix}1&0\\0&p^{\mathbb Z}\end{pmatrix}$} s'identifient
\`a $H^2_{\rm Iw}(\Q,T)$. 
\end{rema}

\subsection{Irr\'eductibilit\'e et admissibilit\'e}\label{resul}

Le but de ce paragraphe est de d\'eduire de tous les r\'esultats pr\'ec\'edents la non nullit\'e, l'irr\'eductibilit\'e (topologique) et l'admissibilit\'e de $\Pi(V)$. 

\begin{coro} \label{cestnonnul}
L'espace de Banach $\Pi(V)$ est non nul.
\end{coro} 
\begin{proof} 
Cela r\'esulte du th\'eor\`eme \ref{phigammabanach} et du corollaire
\ref{iwnonzero2}. 
\end{proof} 
 
Le corollaire \ref{cestnonnul} \'etait conjectur\'e (via la proposition \ref{existencereseau}) et d\'emontr\'e pour $k\leq 2p$ si $p\ne 2$ et $k<4$ si $p=2$ dans \cite[\S3.3]{Br1} par un calcul explicite de r\'eseaux.
 
\begin{coro}\label{irreductible} 
Le $\G$-Banach unitaire $\Pi(V)$ est topologiquement irr\'educti\-ble. 
\end{coro} 
\begin{proof} 
Cela r\'esulte du th\'eor\`eme \ref{phigammabanach} et de la proposition
\ref{actgirred}.
\end{proof} 
 
La proposition \ref{actgirred} montre que $\Pi(V)$ est en fait topologiquement irr\'eductible comme $\B$-repr\'esentation.
 
\begin{coro}\label{admissible} 
Le $\G$-Banach unitaire $\Pi(V)$ est admissible. 
\end{coro} 
\begin{proof} 
On ignore si le $\O$-module compact $\varprojlim_{\psi}D^{\sharp}(T)$ est stable par $\K$ dans $\Pi(V)^*$ (via le th\'eor\`eme \ref{phigammabanach}) mais on peut le remplacer par le $\O$-r\'eseau de $\Pi(V)^*$ :
$${\mathcal M}\=\cap_{g\in \K}g(\varprojlim_{\psi}D^{\sharp}(T))\subset \varprojlim_{\psi}D^{\sharp}(T)$$
qui est un sous-$\O[[X]]$-module compact stable par $\G$ dans $\Pi(V)^*$ (on v\'erifie qu'il est stable par $B(\Q)$ en utilisant la d\'ecomposition d'Iwasawa de $\G$). Le $\O$-module ${\mathcal M}$ poss\`ede alors {\it deux} structures naturelles de $\O[[X]]$-modules : l'une est celle d\'ej\`a d\'efinie et l'autre est :
$$(\lambda,v)\in \O[[X]]\times {\mathcal M}\mapsto
\begin{pmatrix}0&1\\1&0\end{pmatrix}\lambda\begin{pmatrix}0&1\\1&0\end{pmatrix}v.$$
La premi\`ere structure est telle que la multiplication par $(1+X)^{\Z}$ correspond \`a l'action de {\tiny$\begin{pmatrix}1&\Z\\0&1\end{pmatrix}$} et la deuxi\`eme est telle que la multiplication par $(1+X)^{\Z}$ correspond \`a l'action de {\tiny$\begin{pmatrix}1&0\\\Z&1\end{pmatrix}$}. Soit ${\rm pr}:{\mathcal M}\rightarrow D^{\sharp}(T)$ la projection sur la premi\`ere composante et $M\={\rm pr}({\mathcal M})$ : $M$ est un sous-$\O[[X]]$-module (de type fini) de $D^{\sharp}(T)$. Posons ${\mathcal N}\={\rm Ker}({\rm pr})\subsetneq {\mathcal M}$. L'application :
\begin{equation}\label{tordu}
{\mathcal N}\rightarrow M,\ v\mapsto {\rm pr}\left({\tiny \begin{pmatrix}0&1\\1&0\end{pmatrix}}v\right)
\end{equation}
est injective : si $v$ a pour image $0$, sa distribution associ\'ee $\mu_{\alpha}\in B(\alpha)^*\simeq {\mathcal C}^{{\rm val}(\alpha)}(\Z,L)^*\oplus {\mathcal C}^{{\rm val}(\alpha)}(\Z,L)^*$ par (\ref{versbanach}) et (\ref{f1f2}) est nulle sur les deux copies de ${\mathcal C}^{{\rm val}(\alpha)}(\Z,L)$, donc est nulle dans $\Pi(V)^*$. En pensant encore en termes de distributions, on voit que ${\mathcal N}$ est un $\O[[X]]$-module pour la premi\`ere structure mais seulement un $\varphi(\O[[X]])$-module pour la deuxi\`eme structure. De plus, pour cette deuxi\`eme structure, l'injection (\ref{tordu}) est $\varphi(\O[[X]])$-lin\'eaire. Comme $M$ est de type fini sur $\O[[X]]$, donc sur $\varphi(\O[[X]])$, on obtient que le $\varphi(\O[[X]])$-module ${\mathcal N}$ pour la deuxi\`eme action de $\varphi(\O[[X]])$ est de type fini. Fixons maintenant des \'el\'ements $(e_1,\cdots,e_m)\in {\mathcal M}$ (resp. $(f_1,\cdots,f_n)\in {\mathcal N}$) tels que les ${\rm pr}(e_i)$ (resp. les $f_i$) engendrent $M$ sur $\O[[X]]$ (resp. $\mathcal N$ sur $\varphi(\O[[X]])$). Soit $v\in {\mathcal M}$. Il existe $\lambda_1,\cdots,\lambda_m$ dans $\O[[X]]$ tels que $v-\sum \lambda_ie_i\in {\mathcal N}$ et il existe $\mu_1,\cdots,\mu_n$ dans $\varphi(\O[[X]])$ tels que $v-\sum \lambda_ie_i=\sum {\tiny\begin{pmatrix}0&1\\1&0\end{pmatrix}}\mu_i{\tiny\begin{pmatrix}0&1\\1&0\end{pmatrix}}f_i$. Comme les $\lambda_i$ correspondent \`a l'action d'\'el\'ements de l'alg\`ebre de groupe de {\tiny$\begin{pmatrix}1&\Z\\0&1\end{pmatrix}$} et les ${\tiny\begin{pmatrix}0&1\\1&0\end{pmatrix}}
\mu_i{\tiny\begin{pmatrix}0&1\\1&0\end{pmatrix}}$ \`a l'action d'\'el\'ements de l'alg\`ebre de groupe de {\tiny$\begin{pmatrix}1&0\\p\Z&1\end{pmatrix}$}, on voit que $\mathcal M$ est {\it a fortiori} de type fini sur l'alg\`ebre de groupe de $\K$, d'o\`u l'admissibilit\'e. 
\end{proof} 

Les corollaires \ref{irreductible} et \ref{admissible} \'etaient conjectur\'es et d\'emontr\'es par un argument de r\'eduction modulo $p$ pour $k\leq 2p$ (et $k<4$ si $p=2$) dans \cite[\S1.3]{Br2} avec l'hypoth\`ese suppl\'ementaire ${\rm val}(\alpha + \beta)\ne 1$ pour le premier.

On peut d\'eduire des r\'esultats pr\'ec\'edents deux autres corollaires, l'un sur les r\'eseaux dans $\pi(\alpha)$ et $\pi(\beta)$, l'autre sur les vecteurs localement analytiques dans $\Pi(V)$.

\begin{coro}\label{reseaurigolo}
Supposons $\alpha\ne p\beta$, alors $\pi(\alpha)$ (resp. $\pi(\beta)$) poss\`ede des $\O$-r\'eseaux stables par $\G$ et tous les $\O$-r\'eseaux stables par $\G$ dans $\pi(\alpha)$ (resp. $\pi(\beta)$) sont commensurables entre eux. Supposons $\alpha= p\beta$, alors on a le m\^eme r\'esultat pour $\pi(\alpha)$.
\end{coro}
\begin{proof}
L'existence de tels $\O$-r\'eseaux r\'esulte du corollaire \ref{cestnonnul} et de la proposition \ref{existencereseau}. Pour montrer qu'ils sont tous commensurables entre eux, il est \'equivalent de montrer qu'ils sont tous commensurables aux $\O$-r\'eseaux de type fini sur $\O[\G]$. Le $\O$-dual d'un $\O$-r\'eseau stable par $\G$ est toujours contenu dans le $\O$-dual d'un $\O$-r\'eseau de type fini sur $\O[\G]$. Par le th\'eor\`eme \ref{complete} et le corollaire \ref{admissible}, ce dernier dual est de type fini sur l'alg\`ebre de groupe compl\'et\'ee de $\K$. Comme c'est une alg\`ebre noeth\'erienne, il en est de m\^eme du premier dual. Cela entra\^\i ne que le compl\'et\'e de $\pi(\alpha)$ (ou $\pi(\beta)$ si $\alpha\ne p\beta$) par rapport \`a un $\O$-r\'eseau stable par $\G$ quelconque est aussi admissible, et donc topologiquement isomorphe \`a $\Pi(V)$ par le corollaire \ref{irreductible} et le fait que la cat\'egorie des $\K$-Banach admissibles est ab\'elienne (\cite[\S3]{ST3}). Tous les $\O$-r\'eseaux stables par $\G$ induisent donc des normes \'equivalentes sur $\pi(\alpha)$ (ou $\pi(\beta)$ si $\alpha\ne p\beta$) ce qui ach\`eve la preuve.
\end{proof}

Comme dans \cite[\S7]{ST4}, on note $\Pi(V)_{\rm an}$ le sous-$L$-espace vectoriel de $\Pi(V)$ des vecteurs localement analytiques, i.e. des vecteurs $v\in \Pi(V)$ tels que l'application orbite $\G\rightarrow \Pi(V)$, $g\mapsto g\cdot v$ est localement analytique. Il est muni d'une topologie naturelle d'espace localement convexe de type compact (cf. \cite[\S7]{ST4}).

Soit :
$${A}(\alpha)\=\left({\rm Ind}_{\B}^{\G}{\rm nr}(\alpha^{-1})\otimes d^{k-2}{\rm nr}(p\beta^{-1})\right)^{\rm an}$$ 
l'induite parabolique localement analytique au sens de \cite{ST1}. On d\'efinit de m\^eme ${A}(\beta)$ en \'echangeant $\alpha$ et $\beta$. On a des injections naturelles continues $\G$-\'equivariantes ${A}(\alpha)\hookrightarrow B(\alpha)$ et ${A}(\beta)\hookrightarrow B(\beta)$. 

\begin{coro}\label{anal}
Supposons $\alpha\ne p\beta$. On a une injection continue $\G$-\'equiva\-riante :
$${\rm A}(\beta)\oplus_{\pi(\beta)}{\rm A}(\alpha)\hookrightarrow \Pi(V)_{\rm an}$$ 
o\`u $\pi(\beta)$ s'envoie dans $A(\alpha)$ via l'entrelacement (\ref{intertw}).
\end{coro} 
\begin{proof} 
Par \cite[\S4]{ST2}, $\pi(\alpha)$ (resp. $\pi(\beta)$) est le seul sous-objet topologiquement irr\'eductible non nul dans $A(\alpha)$ (resp. $A(\beta)$). Par le th\'eor\`eme \ref{complete}, on d\'eduit que les injections ci-dessus induisent encore des injections ${A}(\alpha)\hookrightarrow B(\alpha)/L(\alpha)$ et ${A}(\beta)\hookrightarrow B(\beta)/L(\beta)$. Le r\'esultat d\'ecoule alors du corollaire \ref{entrepadique}.
\end{proof}

Le corollaire \ref{anal} admet une version lorsque $\alpha=p\beta$ que l'on laisse en exercice au lecteur. Terminons avec une conjecture :

\begin{conj}
Supposons $\alpha\ne p\beta$. L'application ${\rm A}(\beta)\oplus_{\pi(\beta)}{\rm A}(\alpha)\hookrightarrow \Pi(V)_{\rm an}$ du corollaire \ref{anal} est un isomorphisme topologique.
\end{conj}


\begin{thebibliography}{WWWWWW}
 
\bibitem[Ber02]{Be1} Berger L., {\it Limites de repr\'esentations
  cristallines}, \`a para\^\i tre \`a Compositio Math, disponible \`a
  l'adresse : \texttt{www.ihes.fr/${}^{\sim}$lberger}. 

\bibitem[Ber04]{Be2} Berger L., {\it \'Equations diff\'erentielles
  $p$-adiques et $(\varphi,N)$-modules filtr\'es}, pr\'epublication
  2004, disponible \`a l'adresse : \texttt{www.ihes.fr/${}^{\sim}$lberger}.

\bibitem[BB04]{BB} Berger L., Breuil C., {\it Towards a $p$-adic
  Langlands programme}, notes d'un cours donn\'e \`a l'\'Ecole d'\'et\'e de
  Hangzhou (ao\^ut 2004), disponibles \`a l'adresse : 
  \texttt{www.ihes.fr/${}^{\sim}$breuil/publications.html}

\bibitem[Bre03a]{Br1} Breuil C., {\it Sur quelques repr\'esentations modulaires et $p$-adiques de ${\rm GL}_2(\Q)$ II}, J. Institut Math. Jussieu 2, 2003, 23-58.
 
\bibitem[Bre03b]{Br2} Breuil C., {\it Invariant $\mathcal L$ et s\'erie sp\'eciale $p$-adique}, \`a para\^\i tre aux Ann. Scient. E.N.S.
 
\bibitem[Bre03c]{Br3} Breuil C., {\it S\'erie sp\'eciale $p$-adique et cohomologie \'etale compl\'et\'ee}, pr\'epublication 2003, disponible \`a l'adresse : \texttt{www.ihes.fr/${}^{\sim}$breuil/publications.html}.

\bibitem[BM04]{BM} Breuil C., M\'ezard A., {\it En pr\'eparation}.

\bibitem[Bum98]{Bu} Bump D., {\it Automorphic forms and representations}, Cambridge Studies in Advanced Math. 55, Cambridge University Press, 1998.
  
\bibitem[CC99]{CC99} Cherbonnier F., Colmez P., {\it Th{\'e}orie d'Iwasawa des repr{\'e}sentations $p$-adiques d'un corps local}, J. Amer. Math. Soc. 12, 1999, 241--268.

\bibitem[Col03]{Co1} Colmez P., {\it La conjecture de Birch et Swinnerton-Dyer $p$-adique}, S\'eminaire Bourbaki 919, juin 2003.

\bibitem[Col04a]{Co2} Colmez P., {\it Une correspondance de Langlands locale $p$-adique pour les repr\'esentations semi-stables de dimension $2$}, pr\'epublication 2004.

\bibitem[CF00]{CF} Colmez P., Fontaine J.-M., {\it Construction des repr\'esentations $p$-adiques semi-stables}, Inv. Math. 140, 2000, 1--43.

\bibitem[Eme04]{Em} Emerton M., {\it $p$-adic $L$-functions and unitary completions of representations of $p$-adic reductive groups}, pr\'epublication 2004.

\bibitem[Fon90]{F90} Fontaine J-M., {\it Repr{\'e}sentations $p$-adiques des corps locaux I}, The Grothendieck Festschrift, Vol. II, 249--309, Progr. Math. 87, Birkh{\"a}user Boston, Boston, MA 1990.
 
\bibitem[Per94]{BP94} Perrin-Riou B., {\it Th\'eorie d'Iwasawa des
  repr\'esentations $p$-adiques sur un corps local}, Inv. Math. 115, 1994, 81--161.

\bibitem[Sch01]{Sc} Schneider P., {\it Nonarchimedean Functional Analysis}, Springer-Verlag, 2001.

\bibitem[ST01]{ST2} Schneider P., Teitelbaum J. (with an appendix by D. Prasad), {\it $U(\mathfrak g)$-finite locally analytic representations}, Representation Theory 5, 2001, 111--128.
 
\bibitem[ST02a]{ST1} Schneider P., Teitelbaum J., {\it Locally analytic
  distributions and $p$-adic representation theory, with applications
  to ${\rm GL}_{2}$}, J. Amer. Math. Soc. 15, 2002, 443-468.

\bibitem[ST02b]{ST3} Schneider P., Teitelbaum J., {\it Banach space representations and Iwasawa theory}, Israel J. Math. 127, 2002, 359--380.
 
\bibitem[ST03]{ST4} Schneider P., Teitelbaum J., {\it Algebras of $p$-adic distributions and admissible representations}, Inv. Math. 153, 2003, 145--196.

\bibitem[Wa96]{W96} Wach N., {\it Repr{\'e}sentations $p$-adiques
  potentiellement cristallines}, Bull. Soc. Math. France 124, 1996, 375--400.
 
\end{thebibliography}
\end{document}